\documentclass[11pt,twoside,spanish]{article}

\usepackage{
babel,epsfig,graphicx,psfrag,amsfonts}

\newcommand{\dist}{\mbox{$\,$dist$\,$}}

\newcommand{\Z}{\mbox{$Z\!\!\!Z$}}

\begin{document}

\vskip .5cm

\newtheorem{lemma}{Lema}[section]
\newtheorem{theorem}[lemma]{Teorema}
\newtheorem{proposition}[lemma]{Proposici\'{o}n}

\newtheorem{teorema}{Teorema}

\newtheorem{claim}[lemma]{Afirmaci\'{o}n}
\newtheorem{corollary}[lemma]{Corolario}

\newtheorem{corolario}{Corolario}

\newtheorem{definition}[lemma]{Definici\'{o}n}

\newtheorem{remark}[lemma]{Observaci\'{o}n}

.

\begin{center}

 \Huge {Particiones de Markov para difeomorfismos de Anosov}

\vspace{3cm}
\huge
 {Eleonora Catsigeras \footnote{Profesora Agregada del Instituto de Matem\'{a}tica y Estad\'{\i}stica \lq\lq Prof. Rafael Laguardia"(IMERL), Facultad de Ingenier\'{\i}a, Universidad de la Rep\'{u}blica.}}

\end{center}

\normalsize

\pagebreak

 .

 \pagebreak

\begin{center}
{\Large Particiones de Markov para difeomorfismos de Anosov - }
Eleonora Catsigeras
\end{center}

\vspace{.5cm}

\begin{center}

{ \Large \'{I}NDICE}

\vspace{1cm}

\end{center}

0. RESUMEN

\vspace{.6cm}

 1. DEFINICIONES Y RESULTADOS PREVIOS

    \begin{itemize}
    \item[1.1] Difeomorfismos de Anosov
    \item[1.2] Expansividad
    \item[1.3] Estabilidad topol\'{o}gica
    \item[1.4] Conjunto estable e inestable
    \item[1.5] Variedades invariantes
    \item[1.6] Intersecci\'{o}n de variedades invariantes
    \item[1.7] Forma local del producto
    \end{itemize}

    2. RECT\'{A}NGULOS Y PARTICIONES DE MARKOV

 \begin{itemize}
    \item[2.1] Rect\'{a}ngulo de Markov
    \item[2.2] Borde de un rect\'{a}ngulo
    \item[2.3] Propiedades de los rect\'{a}ngulos de Markov
    \item[2.4] Definici\'{o}n de partici\'{o}n de Markov
    \item[2.5] Borde de la partici\'{o}n
    \item[2.6] Propiedades de las particiones de Markov

    \end{itemize}

    3. SEMICONJUGACI\'{O}N CON EL SHIFT

\begin{itemize}
    \item[3.1] Espacio de las sucesiones y funci\'{o}n shift
    \item[3.2] Conjugaci\'{o}n y semiconjugaci\'{o}n
    \item[3.3] Semiconjugaci\'{o}n entre el difeomorfismo de Anosov y
    el shift.
    \item[3.4] Conjuntos estable e inestable en el espacio de las
    sucesiones.
    \item[3.5] Construcci\'{o}n de un cubrimiento con rect\'{a}ngulos
    \item[3.6] Propiedades del cubrimiento.

    \end{itemize}

\pagebreak

    4. M\'{E}TODO CONSTRUCTIVO PARA LA PARTICI\'{O}N

\begin{itemize}
    \item[4.1] M\'{e}todo constructivo para la partici\'{o}n.
    \item[4.2] Segundo refinamiento.
    \item[4.3] Densidad de un conjunto cubierto por el
    refinamiento.
    \item[4.4] Obtenci\'{o}n de la partici\'{o}n.
    \item[4.5] Densidad de las variedades estable e inestable.

    \end{itemize}

    5. TEOREMA DE SINAI

    \begin{itemize}
    \item[5.1] Enunciado
    \item[5.2] Lema
    \item[5.3] Demostraci\'{o}n

    \end{itemize}

6. DIN\'{A}MICA SIMB\'{O}LICA

 \begin{itemize}
    \item[6.1] Matriz de transici\'{o}n
    \item[6.2] Lema
    \item[6.3] Teorema de semiconjugaci\'{o}n
    \item[6.4] Coclusi\'{o}n

    \end{itemize}

    7. REFERENCIAS BIBLIOGR\'{A}FICAS

\vspace{7cm}

    \pagebreak

\begin{center}
{\Large Particiones de Markov para difeomorfismos de Anosov - }
Eleonora Catsigeras
\end{center}

\vspace{3cm}

     \begin{center}
     {\Large RESUMEN }
     \end{center}

\vspace{.5cm}

 En este libro se define Partici\'{o}n de Markov de una
variedad compacta y riemanniana $M$ para un difeomorfismo $f$ de
$M$ en $M$, cuando este difeomorfismo pertenece a cierta clase
particular llamada "de Anosov". El motivo  es demostrar
la existencia de particiones de Markov. Est\'{a} dirigido a estudiantes y egresados
de nivel de grado universitario en Matem\'{a}tica.

   En la parte 1 se exponen definiciones y teoremas que se asumen
   conocidos, referentes a los difeomorfismos llamados de Anosov.
   Casi todos los resultados expuestos en la secci\'{o}n 1 se enuncian
   sin demostraci\'{o}n porque no son el motivo de esta presentaci\'{o}n.
   Pueden encontrarse en las referencias \cite{2}, \cite{3} y
   \cite{4}.

    En la parte 2 se define partici\'{o}n de Markov de una variedad
    $M$ para un difeomorfismo $f$ en $M$. La definici\'{o}n est\'{a}
    referida solamente a los difeomorfismos de Anosov, aunque es
    aplicable a una clase m\'{a}s general de difeomorfismos en la
    variedad. (V\'{e}ase \cite{5}).

    La demostraci\'{o}n de la existencia de una partici\'{o}n de Markov (Teorema de Sinai
    \cite{6}) se concluye en la parte 5 de este trabajo y se basa
    en la construcci\'{o}n de un cubrimiento adecuado de la variedad
    que se luego refinado apropiadamente. Este m\'{e}todo constructivo
    fue extra\'{\i}do del libro de R. Bowen (\cite{1}).

     Tambi\'{e}n se extrajo de R. Bowen (\cite{1}) la presentaci\'{o}n de
     la din\'{a}mica simb\'{o}lica que se expone en la \'{u}ltima parte de esta monograf\'{\i}a.

\vspace{4cm}

\pagebreak

.
 \pagebreak

\begin{center}
{\Large Particiones de Markov para difeomorfismos de Anosov - }
Eleonora Catsigeras
\end{center}

\section{Definiciones y Resultados Previos}
Se asume que $M$ es una variedad de clase $C^1$ compacta y
Riemanniana y que $f:M \mapsto M$ es un dofeomorfismo de clase
$C^1$ en $M$.

\vspace{.3cm}

{\bf \large 1.1 Difeomorfismo de Anosov}

\begin{definition} \em
$f$ es un \em difeomorfismo de Anosov \em  si existen constantes
$K$ y $\lambda$; $K>0$, $0 < \lambda < 1$; y subespacios $S_x,
U_x$ de $T_x M$, tales que:
\begin{itemize}
\item[i] $S_x, U_x$ var\'{\i}an continuamente con $x$
\item [ii] $S_x \oplus U_x = T_xM$
\item[iii] $S_x, U_x$ son invariantes con $f$, es decir:
$f'_xS_x = S_{f(x)}, f'_xU_x = U_{f(x)} \; \forall x \in M$
\item[iv] $\|(f^n)'_x s_x\| \leq K \lambda ^n \|s_x\| \; \forall
s_x \in S_x, \; \forall n \geq 0, \; \forall x \in M$
\item[v] $\|(f^n)'_x u_x\| \leq K \lambda ^{-n} \|u_x\| \; \forall
u_x \in U_x, \; \forall n \leq 0, \; \forall x \in M$
\end{itemize}
\end{definition}

\begin{remark} \em
Las dimensiones de $S_x$ y $U_x$ son constantes en las componentes
conexas de $M$ debido a la condici\'{o}n i. de la definici\'{o}n anterior.
El fibrado tangente $TM$ es la suma directa de los dos subfibrados
$S$ y $U$ invariantes con $f$, llamados \lq \lq subfibrado estable" e
\lq \lq inestable" respectivamente. Si $(x,s_x) \in S$, entonces su norma
decrece (m\'{a}s que exponencialmente con tasa $+\log \lambda <0$)
cuando $n \rightarrow + \infty$. Si $(x,u_x) \in U$, entonces su
norma crece (m\'{a}s que exponencialmente con tasa $-\log \lambda >0$)
cuando $n \rightarrow + \infty$, pues sustituyendo en v. $m = -n,
y = f^{-m}(x), u_y = (f^{-m})'_xu_x$ resulta:
$$\|(f^m)'_y u_y\| \geq \frac{1}{K} \lambda ^{-m}\|u_y\|\; \; \forall u_y \in U_y, \; \forall
m \geq 0, \; \forall y \in M$$ Si $f$ es un difeomorfismo de
Anosov, tambi\'{e}n lo es $f^{-1}$, y el subfibrado estable para
$f^{-1}$ es el inestable para $f$ y viceversa, como se observa de
la definici\'{o}n de difeomorfismo de Anosov.
\end{remark}

{\bf \large 1.2 Expansividad}

\begin{definition} \label{121}
\em Un difeomorfismo $f: M \mapsto M$ es \em expansivo \em si
existe una constante $\rho >0$, llamada \em constante de
expansividad, \em tal que
$$\dist (f^nx, f^ny)\leq \rho \; \forall n \in \Z \;\; \mbox {si y solo si }
\;\; x= y $$ Una sucesi\'{o}n bi-infinita $\{y_n\}_{n \in \Z}, \; y_n
\in M $, se dice que \em  $\epsilon-$acompa\~{n}a \em a otra
$\{x_n\}_{n \in \Z}, \; x_n \in M$, si $\dist (x_n, y_n) \leq
\epsilon \; \forall n \in \Z$. La expansividad de un difeomorfismo
significa que para cierto $\rho >0 $ suficientemente peque\~{n}o, dos
\'{o}rbitas diferentes nunca se $\rho-$acompa\~{n}an.
\end{definition}
Un propiedad conocida es la siguiente:
\begin{proposition}
Todo difeomorfismo de Anosov es expansivo.
\end{proposition}
{\bf \large 1.3 Estabilidad topol\'{o}gica.}
\begin{definition}
\em Un difeomorfismo $f: M \mapsto M$ es \em topol\'{o}gicamente
estable \em si dado $\epsilon > 0$ existe un $C^0$ entorno ${\cal
V}$ de $f$ tal que para todo $g \in {\cal V}$ existe una
semiconjugaci\'{o}n $h: M \mapsto M$ entre $g$ y $f $ (i.e. h es
continua, sobreyectiva y cumple $h \circ g = f \circ h$) tal que
$\dist (h(x),x) < \epsilon \;\; \forall x \in M$.
\end{definition}
\begin{remark}
Si $f$ es topol\'{o}gicamente estable, si ${\cal V}$ es un $C^0$
entorno de $f$ como en la definici\'{o}n anterior y si $g \in {cal V}$
entonces
$$h (g^n(x))= f^n(h(x)) \;\; \forall x \in M \;\; \forall n \in \Z$$
La estabilidad topol\'{o}gica de $f$ significa que si $g$ est\'{a}
suficientemente pr\'{o}xima de $f$ (en la topolog\'{\i}a $C^0$) entonces
las \'{o}rbitas de $g$ est\'{a}n $\epsilon$-acompa\~{n}adas por las de $f$ y
$\epsilon$ acompa\~{n}an a todas las \'{o}rbitas de $f$ (ya que la
transformaci\'{o}n $h$ es sobreyectiva.)
\end{remark}
\begin{theorem} \label{133}
\em {\bf (Pugh) } \em Un difeomorfismo $f: M \mapsto M$ es
topol\'{o}gicamente estable si y s\'{o}lo si dado $\epsilon >0$ existe
$\delta >0$ tal que toda sucesi\'{o}n bi-infinita de puntos $\{y_n\}{n
\in \Z}, \; y_n \in M$ que cumple $\dist (y_{n+1}, f(y_n))< \delta
\; \forall n \in \Z$ est\'{a} $\epsilon$-acompa\~{n}ada por una \'{o}rbita de
$f$
\end{theorem}
\begin{definition}
\em Sea $f: M \mapsto M$ invertible. Una sucesi\'{o}n bi-infinita
$\{y_n\}_{n \in \Z}$ de puntos de $M$ se llama \em
$\delta$-pesudo-\'{o}rbita de $f$ \em si $\dist (f(y_n), y_{n+1})<
\delta \; \forall n \in \Z$.
\end{definition}
Se concluye que un difeomorfismo $f: M \mapsto M$ es
topol\'{o}gicamente estable si y solo si dado $\epsilon >0$ existe
$\delta >0$ tal que  toda $\delta$-pseudo-\'{o}rbita est\'{a}
$\epsilon$-acompa\~{n}ada.
\begin{theorem} \label{135}
Los difeomorfismos de Anosov son topol\'{o}gicamente estables.
\end{theorem}
{\bf \large 1.4 Conjuntos estable e inestable.}
\begin{definition} \label{141}
\em Sea $f: M \mapsto M$ un difeomorfismo. Se llama \em conjunto
estable \em de $f$ por el punto $x \in M$ a
$$W^s(x)= \{y \in M: \dist(f^n(y), f^n(x)) _{n \rightarrow + \infty} \rightarrow 0\}$$
Se llama \em conjunto inestable \em de $f$ por el punto $x \in M$
a
$$W^u(x)= \{y \in M: \dist(f^n(y), f^n(x)) _{n \rightarrow - \infty} \rightarrow 0\}$$
\end{definition}
\begin{theorem}
Si $f$ es un difeomorfismo de Anosov entonces  \em $$W^s=\{ y \in
M: \limsup_{n \rightarrow + \infty} \frac {1}{n} \log  \dist
(f^n(x), f^n(y)) \leq \log \lambda\}$$
$$W^u=\{y \in M:
\limsup_{n \rightarrow + \infty} \frac {1}{n} \log  \dist
(f^{-n}(x), f^{-n}(y)) \leq \log \lambda\}$$  \em  donde $\lambda$
es la misma constante $0< \lambda < 1$ de la definici\'{o}n de Anosov.
\end{theorem}
El teorema anterior significa que para los difeomorfismos de
Anosov, dos \'{o}rbitas que en el futuro (o en el pasado) se acercan
de modo que su distancia tienda a cero, entonces lo hacen m\'{a}s que
exponencialmente con tasa $\log \lambda$. En efecto: $\lim sup_{n
\rightarrow + \infty} \frac {1}{n} \log \dist (f^n(x), f^n(y))
\leq \log \lambda$ implica $$\dist (f^n(x),f^n(y))< A e ^{-n
\gamma} \; \; \forall n \geq 0$$ donde $A$ es un n\'{u}mero positivo y
$\gamma$ es un n\'{u}mero real positivo elegido de modo que $\gamma <
-\log \lambda$.

\begin{remark} \em
Dos conjuntos estables distintos son disjuntos pues si $z \in
W^s(x)\cap W^s(x') \neq \emptyset$ entonces, a partir de la
definici\'{o}n de conjunto estable y la propiedad triangular de la
distancia se tiene  $W^s(z) \subset W^s(x) \cap W^s(x')$ y
 $W^s(x) \cup W^s(x') \subset W^s(z)$. Luego $W^s(x) = W^s(x')$.

 Los mismo vale para los conjuntos inestables.
\end{remark}
\begin{proposition} \label{145}
Si $f$ es un difeomorfismo expansivo con constante de expansividad
$\rho$, entonces: \em
$$W^s(x)= \{y \in M: \dist (f^n(x), f^n(y)) \leq \rho \; \forall n \mbox{ suficientemente
grande }\}$$
$$W^u(x)= \{y \in M: \dist (f^{-n}(x), f^{-n}(y)) \leq \rho \; \forall n
\mbox{ suficientemente
grande }\}$$
\end{proposition}

{\em Demostraci\'{o}n: }De la definici\'{o}n de conjunto estable se
obtiene que $$W^s(x) \subset \{y\in M: \dist (f^n(x),f^n(y)) \leq
\rho \; \forall n \mbox{ suficientemente grande}\})$$ Para
demostrar la otra inclusi\'{o}n supongamos por absurdo que existe $y
\in M$ tal que
$$\dist (f^n(x), f^n(y)) \leq \rho \; \forall n \geq N , \; \dist (f^n(x),f^n(y))_{
n \rightarrow + \infty} \not \rightarrow 0$$ Entonces existe una
sucesi\'{o}n $n_k \rightarrow + \infty$ tal que $\dist(f^{n_k}(x),
f^{n_k}(y))\geq \epsilon >0 \; \forall k$. Podemos elegir $n_k$ de
modo que $f^{n_k}(x)$ y $ f^{n_k}(y)$ sean convergentes (por la
compacidad de $M$). Luego si $p \in \Z$ se tiene
$$\dist (f^{n_k+p}(x),f^{n_k+p}(y))\leq \rho \;\; \forall n_k >N -p $$
Cuando $k \rightarrow + \infty$ se tiene $\dist (f^p(x_0),
f^p(y_0))\leq \rho \;\forall p \in \Z$ donde $x_0 = \lim
f^{n_k}(x), \; y_0 = \lim f^{n_k}(y)$. Por la expansividad $x_0 =
y_0$, contradiciendo la elecci\'{o}n de $n_k$, pues $$
\dist(f^{n_k}(x), f^{n_k}(y))\geq \epsilon >0 \; \;\;\;  \Box$$

\begin{remark} \label{146} \em
Se observa que la proposici\'{o}n anterior sigue siendo v\'{a}lida si se
sustituye la constante de expansividad $\rho$ por cualquier otra
constante $\epsilon
>0, \epsilon \leq \rho$. Entonces
$$W^s(x) = \bigcup _{n \in \Z} \{y \in M: \dist (f^n(x), f^n(y))\leq \epsilon \;
\forall n \geq N\}=$$ $$ = \bigcup _{n \geq k} \{y \in M: \dist
(f^n(x), f^n(y))\leq \epsilon \; \forall n \geq N\} \;\; \;\forall
k \in \Z $$
\end{remark}
\begin{definition} \em
Se llama \em $\epsilon-$conjunto estable \em de $f$ por el punto
$x \in M$ al conjunto
$$W^s_\epsilon (x)= \{y \in M: \dist (f^n(x), f^n(y))\leq \epsilon \; \forall n \geq 0\}$$

Se llama \em $\epsilon-$conjunto inestable \em de $f$ por el punto
$x \in M$ al conjunto
$$W^u_\epsilon (x)= \{y \in M: \dist (f^n(x), f^n(y))\leq \epsilon \; \forall n \leq 0\}$$
\end{definition}
\begin{remark} \em \label{148}
Si $f$ es un difeomorfismo expansivo es f\'{a}cil ver, a partir de la
definici\'{o}n anterior y de \ref{145} y \ref{146} que:
\begin{itemize}
\item[i] $y \in W^s_{\epsilon } (x) \; \Leftrightarrow \; x \in W^s_{\epsilon}
(y); \;\;\;\; y \in W^u_{\epsilon } (x) \; \Leftrightarrow \; x
\in W^u_{\epsilon} (y)$.
\item[ii] Si $\epsilon \leq \rho$ entonces $ W^s_\epsilon (x) \cap
W^u_\epsilon (x) = \{x\}$.
\item[iii] $f(W^s_\epsilon(x)) = \{z \in M: \dist (f^{n}(x), f^{n-1}(z))\leq \epsilon \;
\forall n \geq 0\} = $ $$ = \{z \in M: \dist (f^{n-1}(f(x)),
f^{n-1}(z))\leq \epsilon \; \forall n \geq 1\}\cap \{z \in M:
\dist (x, f^{-1}(z))\leq \epsilon \} = $$ Luego $$
f(W^s_\epsilon(x)) = W^s_{\epsilon }(f(x)) \cap f(\overline
B_{\epsilon }(x))$$ An\'{a}logamente $$ W^u_\epsilon(f(x)) =
f(W^u_{\epsilon }(x)) \cap \overline B_{\epsilon }(f(x))$$
\item[iv] $f(W^s_{\epsilon}(x)) \subset W^s_{\epsilon }(f(x))$ y
  $f(W^u_{\epsilon}(x)) \supset
W^u_{\epsilon }(f(x))$
\item[v] Si $\epsilon \leq \rho$ (v\'{e}ase
\ref{146}) entonces: $$W^s (x) = \bigcup _{N \geq 0}
f^{-N}(W^s_{\epsilon} (f^N(x))) = \bigcup _{N \in \Z}
f^{-N}(W^s_{\epsilon } (f^N(x)))$$  En particular
$$W^s_{\epsilon }(x) \subset W^s(x),
 \;\;\; W^u_{\epsilon }(x) \subset W^u(x)$$
\end{itemize}
\end{remark}
{\bf \large 1.5 Variedades invariantes.}

\vspace{.5cm} El teorema que sigue justifica el nombre de \em
variedad invariante estable \em (respectivamente \em inestable\em
) que recibe el conjunto estable (respectivamente inestable)
cuando $f$ es un difeomorfismo de Anosov. Se enuncia sin
demostraci\'{o}n, la cual puede encontrarse en la referencia \cite{3}.
\begin{theorem} \label{151}
Sea $f: M \mapsto M$ un difeomorfismo de Anosov. Entonces los
conjuntos $W^s(x)$ y $W^u (x)$ son $C^1$ variedades inmersas en
$M$, que pasan por $x$, tangentes en $x$ a los subespacios $S_x$ y
$U_x$ respectivamente.
\end{theorem}
Se observa de la definici\'{o}n \ref{141} que la partici\'{o}n de $M $ en
las variedades estables e inestables es invariante por $f$; m\'{a}s
precisamente: $f(W^s(x)) = W^s(f(x)); \; f(W^u(x)) = W^u(f(x))$
\begin{remark} \em
Cuando $f$ es un difeomorfismo de Anosov, entonces $W^s(x)$ y
$W^u(x)$ son variedades inmersas en $M$ seg\'{u}n afirma el teorema
\ref{151}. Sin embargo no son necesariamente subvariedades de $M$:
la inclusi\'{o}n es continua pero no necesariamente un homeomorfismo
sobre su imagen (la topolog\'{\i}a en $W^s(x)$ podr\'{\i}a ser estrcitamente
m\'{a}s fina que la inducida por $M$).
\end{remark}
\begin{theorem}
Sea $f: M \mapsto M$ un difeomorfismo de Anosov. Si $\epsilon >0$
es suficientemente peque\~{n}o entonces:
\begin{itemize}
\item [1] Para todo $N \in \Z$ el conjunto
$f^N(W^s_{\epsilon }(f^{-N}(x)))$ es un entorno de $x$ en la
variedad $W^s(x)$ homeomorfo a una bola en el subespacio $S_x$
\item [2] La topolog\'{\i}a en $f^N(W^s_{\epsilon }(f^{-N}(x)))$ como
entorno en la variedad $W^s(x)$ es la misma que la inducida en \'{e}l
por la topolog\'{\i}a de $M$
\end{itemize}

\end{theorem}

La demostraci\'{o}n se encuentra en l referencia bibliogr\'{a}fica [3]. Es
parte de la demostraci\'{o}n del teorema  \ref{151} de existencia
variedades invariantes.

 En particular
$W^s_{\epsilon } (x)$ es un entorno de $x$ en $W^s(x)$ y la
topolog\'{\i}a en $W^s_{\epsilon }(x)$ como subconjunto de $W^s(x)$ es
la misma que la inducida por la topolog\'{\i}a de $M$ como subconjunto
de $M$.

Aplicando el teorema anterior a $f^{-1}$ en lugar de $f$ y
observando que por definici\'{o}n las variedades estables de $f^{-1}$
son las inestables de $f$, se obtiene que $W^u_{\epsilon } (x)$ es
un entorno de $x$ en $W^u(x)$ y su topolog\'{\i}a como subconjunto de
$W^u(x)$ es la misma que la inducida por la topolog\'{\i}a de $M$.

En virtud de la primera parte del teorema anterior, los $epsilon-$
conjuntos estable e inestable se llaman $epsilon-$variedades
estable e inestable, y son subvariedades de $M$ (variedades
encajadas, con la topolog\'{\i}a inducida por la de $M$).

De la segunda parte del teorema anterior se desprende que:
\begin{itemize}
\item Si $U$ es un abierto de $W^s_{\epsilon }(x)$ (como variedad estable), entonces
existe $B$ abierto en $M$ tal que $U = B \cap W^s_{\epsilon }(x)$
\item Si $x_n, x \in W^s_{\epsilon }(x)$ y si $\dist (x_n, x) \rightarrow
0$ en $M$, entonces $x_n \rightarrow x$ en $W^s_{\epsilon }(x)$ (y
rec\'{\i}procamente).

\end{itemize}
{\bf \large 1.6 Intersecci\'{o}n de variedades invariantes.}

\vspace{.3cm}

 En la observaci\'{o}n \ref{148} obtuvimos que
$W^s_{\epsilon } (x) \cap W^u{\epsilon } (x) = \{x\}$ si $0
<\epsilon < \rho$. Probaremos a continuaci\'{o}n que cuando $x$ e $y$
est\'{a}n suficientemente pr\'{o}ximos y $\epsilon >0 $ es peque\~{n}o,
entonces $W^s_{\epsilon }\cap W^u{\epsilon } (y) $ consiste en un
\'{u}nico punto que denotaremos como $[x,y]$. En lo que sigue $f: M
\mapsto M$ denota un difeomorfismo de Anosov.

\begin{proposition} \label{161}
Dado $0<\epsilon \leq \rho /2$  existe $\delta >0$ tal que
$W^s_{\epsilon } (x) \cap W^u_{\epsilon }(y)$ consiste en un \'{u}nico
punto, para todos $x,y \in M$ tales que $\dist (x,y) < \delta$.
\end{proposition}
{\em Demostraci\'{o}n: } Por los teoremas \ref{133} y \ref{135} existe
$\delta _1
>0 $ tal que toda $\delta _1$ pseudo-\'{o}rbita de $f$ est\'{a} $\epsilon
/2$ acompa\~{n}ada por una \'{o}rbita de $f$.

Tomemos $delta = \min (\delta _1, \epsilon /2)$ y dos puntos $x,y
\in M$ tales que $\dist (x,y) < \delta$.

La sucesi\'{o}n bi-infinita $\{y_n\}_{n \in \Z}$ definida por $y_n =
f^n(x)$ si $n \geq 0 $, $ y_n = f^n (y_n)$ si $n < 0 $ es una
$\delta _1$ pseudo-\'{o}rbita. Entonces existe $z \in M$  que cumple
$\dist (f^n(z), f^n(x)) \leq \epsilon /2 \; \forall n \geq 0,
\dist (f^n(z), f^n(y)) \leq \epsilon /2 \; \forall n < 0$. Adem\'{a}s
$\dist (z,y) \leq \dist (z,x) +   \dist (x,y) \leq \epsilon /2 +
\delta \leq \epsilon$.

Entonces $z \in W^s_{\epsilon }(x) \cap W^u_{\epsilon }(y)$.

Adem\'{a}s $z$ es \'{u}nico en $W^s_{\epsilon }(x) \cap W^u_{\epsilon
}(y)$ debido a la expansividad de $f$. $\;\; \Box$
\begin{definition}
\em Llamaremos funci\'{o}n \em corchete \em a
$$[\cdot, \cdot]: \{(x,y) \in M^2: \dist (x,y < \delta)\} \mapsto M$$
definida por $$[x,y]= W^s_{\epsilon }(x) \cap W^u_{\epsilon }(y)$$
\end{definition}

Se observa que $$\dist (f^n([x,y]), f^n(x))\leq \epsilon \;
\forall n \geq 0, \;\;\; \dist (f^n([x,y]), f^n(y))\leq \epsilon
\; \forall n \leq 0$$ debido a la definici\'{o}n de la funci\'{o}n
corchete y a la definici\'{o}n de las $\epsilon$-variedades estable e
inestable.
\begin{theorem}
La funci\'{o}n corchete $[\cdot, \cdot]$ es continua.
\end{theorem}
{\em Demostraci\'{o}n: } Sea $(x_n, y_n ) \rightarrow (x,y)$ en $M^2$
tales que $\dist (x,y) < \delta$. Como $M$ es compacta puede
elegirse $(x_n, y_n)$ de modo que  $[x_n,y_n]$ sea convergente en
$M$. Sea $z_n= [x_n,y_n] \rightarrow z \in M$. Basta demostrar que
$z = [x,y]$.

Como $z_n \in W^s_{\epsilon }(x_n)$ entonces $\dist (f^p(z_n),
f^p(x_n))\leq \epsilon \; \forall p \geq 0$. Dejando fijo $p$ y
haciendo $n \rightarrow \infty$, por la continuidad de $f$ se
tiene que $\dist (f^p(z), f^p(x))\leq \epsilon \; \forall p \geq
0$, de donde $z \in W^s_{\epsilon }(x)$.

An\'{a}logamente se obtiene $z \in W^u_{\epsilon }(y)$, de donde $z
\in W^s_{\epsilon } (x) \cap W^u_{\epsilon }(y) = [x,y] \;\; \Box$
\vspace{.4cm}

{\bf \large 1.7 Forma local del producto}
\begin{theorem} \label{171}
Sea $f: M \mapsto M$ un difeomorfismo de Anosov. Existe una
constante $\overline \epsilon >0 $ tal que para todo $x \in M$ el
producto $W^u_{\overline \epsilon }(x) \times W^s _{\overline
\epsilon }(x)$ es homeomorfo a un entorno de $x $ en $M$.
\end{theorem}
{\em Demostraci\'{o}n: } Elijamos $0<\epsilon < \rho /2$ (donde $\rho$
es la constante de expansividad de $f$). Sea $\delta >0$ elegido
seg\'{u}n el teorema \ref{161} y sea $0 < \overline \epsilon < \min
(\delta /2, \epsilon )$. Resulta:
$$W^u_{\overline \epsilon } (x) \times W^s_{\overline \epsilon} (x)
\subset \{(z,y) \in M^2: \dist (z,y)< \delta \} $$
 Puede aplicarse la funci\'{o}n corchete a puntos en $W^u_{\overline \epsilon}
  (x) \times W^s_{\overline \epsilon} (x)
 \subset M^2$.

  Sea $\varphi = \left . [\cdot, \cdot]\right |_{
 W^u_{\overline \epsilon  } (x) \times W^s_{\overline \epsilon}
 (x)}$. La aplicaci\'{o}n $\varphi$ es continua porque es la
 restricci\'{o}n de una funci\'{o}n continua. Es inyectiva pues si $[z,z'] = [\overline z,
 \overline z ']$ donde $z, \overline z \in W^u_{\overline \epsilon}(x), \;
 z', \overline z' \in W^s_{\overline \epsilon}(x)$ entonces
 $$\dist (f^n(z), f^n(\overline z))\leq 2 \overline \epsilon < 2 \epsilon \;
 \forall n \geq 0$$
 $$\dist (f^n(z), f^n(\overline z))\leq 2 \overline \epsilon < 2 \epsilon \;
 \forall n \leq 0$$
 Luego por la expansividad de $f$ se tiene $z = \overline z$.
 An\'{a}logamente se obtiene $z' = \overline z'$.

 La $\overline \epsilon$-variedad estable $W^s_{\overline \epsilon}(x)=
 \{y \in M: \dist (f^n(y), f^n(x))\leq \overline \epsilon \; \forall n \geq
 0\}$ es cerrada en $M$ que es compacta, luego es compacta.
 An\'{a}logamente es compacta $W^u_{\overline \epsilon}(x)$. As\'{\i}
 $\varphi$ es continua e inyectiva con dominio compacto
 $W^u_{\overline \epsilon}(x)\times W^s_{\overline
\epsilon}(x)$  a $M$. Como el dominio de $\varphi$ y su codominio
son variedades de la misma dimensi\'{o}n finita,  $\varphi$ es un
 homeomorfismo sobre su imagen.

 En efecto $\varphi (x,x) = x,
\;  \dim W^s_{\overline \epsilon}(x) = \dim S_x, \;   \dim
W^u_{\overline \epsilon}(x) = \dim U_x, \; S_x \oplus U_x = T_x M$
de donde $\dim (W^u_{\overline \epsilon}(x)\times W^s_{\overline
\epsilon}(x))= \dim M$.

Siendo $W^u_{\overline \epsilon}(x)\times W^s_{\overline
\epsilon}(x)$ un entorno de $(x,x)$ en $M^2$, su imagen homeomorfa
es un entorno de $x \in M. \;\; \Box$

\vspace{4cm}

\pagebreak

\begin{center}
{\Large Particiones de Markov para difeomorfismos de Anosov - }
Eleonora Catsigeras
\end{center}

\section{Rect\'{a}ngulos y particiones de Markov}

Sea $f$ un difeomorfismo de Anosov en una variedad compacta y
Riemanniana $M$. Se definir\'{a} \em Partici\'{o}n de Markov \em para $f$.
Es un cubrimiento finito de $M$ por cierta clase de cerrados
llamados \em rect\'{a}ngulos, \em con interiores dos a dos disjuntos,
y que cumplen condiciones que los vinculan a la din\'{a}mica de $f$,
es decir, al espacio de \'{o}rbitas de $f$.

 Comenzaremos definiendo \em rect\'{a}ngulo \em y demostrando algunas
 propiedades que ser\'{a}n utilizadas m\'{a}s adelante.

\vspace{.3cm}

 {\bf \large 2.1 Definici\'{o}n de rect\'{a}ngulo}

\vspace{.3cm}

  Sea $\epsilon >0$ tal que $0 < \epsilon \rho /4$ (donde $\rho $ denota la constante de
   expansividad de $f$, definida en \ref{121}. y sea $\delta >0$
   elegido como en \ref{161}.

\begin{definition}
\em Un subconjunto $R$ no vac\'{\i}o de $M$ se llama \em rect\'{a}ngulo \em
para  $f $ si tiene di\'{a}metro menor que $\delta$ y adem\'{a}s $[x,y]\in
R\; \forall x,y \in R$. Es decir: si $x,y \in R$ entonces existe
un \'{u}nico $z \in W^s_{\epsilon }(x) \cap W^s_{\epsilon }(y) =
[x,y]$ y adem\'{a}s $z \in R$.
 \end{definition}
\begin{definition}
\em Un rect\'{a}ngulo $R$ es propio si $R = \overline{ \mbox{int}R }$.
\end{definition}

{\bf Ejemplos: }
\begin{itemize}
 \item  [i) ] A partir de la definici\'{o}n obs\'{e}rvese
  que si $\overline \epsilon >0$ se
 elige suficientemente peque\~{n}o entonces es un rect\'{a}ngulo
 el entorno $V$ de $x$ en $M
 $ que tiene
 forma local del producto seg\'{u}n el teorema \ref{171} (es decir $V$ es homeomorfo a
 $W^u_{\overline \epsilon }(x) \times W^s_{\overline \epsilon
 }(x)$).
 \item  [ii) ] Si $x$ e $y$ son dos puntos pr\'{o}ximos en $M$ y si $U$
 y $V$ son entornos de $x$ e $y$ respectivamente, como en el
 ejemplo anterior, entonces $U \cup V \cup [U,V] \cup [V,U]$ es un
 rect\'{a}ngulo.
 \item  [iii) ] Si $U$ y $V$ son dos rect\'{a}ngulos no disjuntos,
 entonces $U \cap V$ es un rect\'{a}ngulo.

\end{itemize}

\begin{definition}
\em Sea $R$ un rect\'{a}ngulo y $x \in R$. Se llama
$\epsilon$-variedad estable de $x$ en $R$ a
$$W^s(x,R)= W^s_{\epsilon }(x) \cap R$$
An\'{a}logamente:
$$W^u(x,R)= W^u_{\epsilon }(x) \cap R$$
\end{definition}

\begin{proposition} \label{215}
Sean $x,y \in R$, $R$ rect\'{a}ngulo. Entonces $y \in W^s(x,R)$ si y
solo si $W^s(x,R) = W^s(y,R)$.
\end{proposition}

{\bf Demostraci\'{o}n: } Como $y \in W^s(y,R)$ es inmediato que
$W^s(x,R ) \subset W^s(y, R)$ implica $y \in W^s(x,R)$.

Para el rec\'{\i}proco alcanza probar que $y \in W^s(x,R)$ implica
$W^s(x,R) \subset W^s(y,R)$ (pues por simetr\'{\i}a $y \in W^s(x,R)$ si
y  solo si $x \in W^s(y,R)$).

Probemos entonces que $y,z \in W^s(x, R)$ implica $z \in
W^s(y,R)$:

Sea $w = [z,y]= W^s_{\epsilon }(z) \cap W^u_{\epsilon } (y)$. Como
$z \in W^s_{\epsilon } (x)$. Entonces $w = W^s_{2 \epsilon }(x )
\cap W^u_{\epsilon} (y)$.

Adem\'{a}s $y \in W^s_{epsilon } (x)$. Entonces $w = W^s_{3 \epsilon
}(y ) \cap W^u_{\epsilon} (y)$. Siendo $3 \epsilon > \rho$, donde
$\rho$ es la constante de expansividad, se tiene que $w = y$ o
sea:
$$y = [z,x] = W^s_{\epsilon }(z) \cap W^u_{\epsilon } (y)$$
de donde $y \in W^s_{epsilon } (z)$, o lo que es lo mismo $z \in
W^s_{\epsilon } (y)$ como quer\'{\i}amos demostrar. $ \Box $

\vspace{.3cm} {\bf \large 2.2 Borde de un rect\'{a}ngulo }

\begin{definition}
\em Se llama \em Borde Estable \em de un rect\'{a}ngulo $R$ al
conjunto $$\partial ^s R = \{x \in R: x \not \in
\mbox{int}W^u(x,R) \mbox { en } W^u_{\epsilon }(x) \}$$ Se llama
\em Borde Inestable \em de un rect\'{a}ngulo $R$ al conjunto
$$\partial ^u R = \{x \in R: x \not \in \mbox{int}W^s(x,R) \mbox {
en } W^s_{\epsilon }(x) \}$$
\end{definition}

Demostraremos que para los rect\'{a}ngulos $R$ cerrados, el borde
topol\'{o}gico de $R$ es $\partial ^s R \cup \partial ^u R$. Adem\'{a}s
para justificar el nombre de borde estable, demostraremos que
$\partial ^s(R)$ est\'{a} formado por la uni\'{o}n de
$\epsilon$-variedades estables en $R$:

\begin{proposition}
$y \in \partial ^s R$ implica $W^s(y,R)\subset \partial ^sR$.

$y \in \partial ^u R$ implica $W^u(y,R)\subset \partial ^uR$
\end{proposition}

{\bf Demostraci\'{o}n: } Por absurdo sea $x \in W^s(y, R)\setminus
\partial ^sR$. Entonces $x \in \mbox {int} W^u (x,R)$ en $W^u
_{\epsilon } (x)$, o sea existe un entorno $V$ de $x$ en
$W^u_{\epsilon }(x)$ contenido en $R$.

Sea $\varphi (z)= [z,x]= W^s_{\epsilon }(z) \cap W^u _{\epsilon}
(x)$ definido para los puntos $z \in W^u {\epsilon }(y)$ que est\'{a}n
a distancia menor que $\delta$ de $x$.

$\varphi$ es continua pues es la restricci\'{o}n de $[\cdot, \cdot]$.

Adem\'{a}s $\varphi (y)= [y,x]= x$ porque $x \in W^s_{\epsilon }(y)$.

Entonces $\varphi ^{-1}(V)$ es un abierto de $W^u_{\epsilon }(y)$
que contiene a $y$. Adem\'{a}s si $z \in \varphi ^{-1} (V)$ entonces
$\varphi (z) = u \in V \subset R$. Luego $[z,x]= u$, de donde $z
\in W^s_{\epsilon } (u)$. Como $z \in W^u _{\epsilon }(y)$ se
obtiene que $z = [u,y]$ con $u,y \in R$. Entonces por definici\'{o}n
de rect\'{a}ngulo $z \in R$. Se tiene as\'{\i} que $\varphi ^{-1} (V)
\subset R$.

Se ha hallado un entorno de $y$ en $W^u_{\epsilon }(y)$ contenido
en $R$. Entonces $y \in \mbox {int} W^u(y,R)$ en $W^u_{\epsilon
}(y)$, o sea $y \not \in \partial ^sR$ contradiciendo la
hip\'{o}tesis. $\Box$

\begin{proposition} \label{223}
Si $R$ es un rect\'{a}ngulo cerrado entonces $\partial R = \partial ^s
R \cup \partial ^u R$.
\end{proposition}

{\bf Demostraci\'{o}n: } Veremos que $\mbox{int}R = R \setminus
(\partial ^sR \cup \partial ^u R)$ Seg\'{u}n la definici\'{o}n de borde
estable e inestable tenemos que
$$R \setminus (\partial ^s R \cup \partial ^u R) =
\{ x \in R : x \in \mbox{int } W^u(x,R) \mbox{ en } W^u_{\epsilon
} (x); \;\; x \in \mbox{int }  W^s(x,R) \mbox{ en } W^s_{\epsilon
}(x)\}$$ Sea $y \in \mbox{int }R$, sea $B$ un entorno de $y$ en
$M$ contenido en $R$. Tenemos que
$$y \in W^u_{\epsilon } (y) \cap B \supset W^u_{\epsilon }(y) \cap R = W^u(y,R)$$
Luego $W^u_{\epsilon }(y) \cap B$ es un entorno de y en
$W^u_{\epsilon }(y)$ que est\'{a} contenido en $W^u(y,R)$, y entonces
$y \not \in \partial ^sR$.

De igual forma tenemos que $y \not \in \partial ^u R$ , con lo
cual deducimos que
$$\mbox{int }R \subset R \setminus (\partial ^sR \cup \partial ^uR)$$
Rec\'{\i}procamente: Si $y \in \mbox{int }W^u(y,R)$ en $W^u_{\epsilon
}(y)$ entonces existe un entorno de $y$ en $W^u_{\epsilon }(y)$
que est\'{a} contenido en $R$. Llamemos $V$ a su intersecci\'{o}n con
$W^u{\overline \epsilon }(y)$, siendo $\overline \epsilon $ >0
elegido como en el Teorema \ref{171}.

$V$ es un entorno de $y$ en $W^u_{\epsilon }(y)$ porque
$W^u{\overline \epsilon }(y)$ lo es. Adem\'{a}s $V \subset R \cap
W^u{\overline \epsilon }(y)$. De igual forma hallemos $U$ entorno
de $y$en $W^s_{\epsilon }(y)$ contenido de $R \cap W^s{\overline
\epsilon }(y) $.

 Por el Teorema \ref{171} $[V,U]$ es un entorno de $y$ en $M$.
 Como $V,U \subset R$, por la definici\'{o}n de rect\'{a}ngulo se deduce
 que $[V,U] \subset R$. Entonces $y \in \mbox{int }R$. Luego
 $R \setminus (\partial ^s R \cup \partial ^u R) \subset \mbox{int
 }R.
 \; \;\;\ \Box$

 \vspace{.3cm}

 {\bf \large 2.3 Propiedades de los rect\'{a}ngulos }

 \vspace{.3cm}

 Algunas propiedades que se demuestran a continuaci\'{o}n
 ser\'{a}n utilizadas en los par\'{a}grafos siguientes:

 \begin{proposition}
 Si $R$ es un rect\'{a}ngulo entonces tambi\'{e}n los son
  $\overline R$ y $\mbox{int }R$ cuando no es vac\'{\i}o.

  Si $R_1$ y $R_2$ son rect\'{a}ngulos tales que $R_1 \cap R_2 \not =
  \emptyset$ entonces $R_1 \cap R_2$ tambi\'{e}n es un rect\'{a}ngulo.
 \end{proposition}
 {\em Demostraci\'{o}n: } Sean $x,y \in \overline R, \; x_n = \rightarrow x, \;
 y_n \rightarrow y, \; x_n, y_n \in R$. Probemos que $[x,y]\in \overline R for
 all x,y \in \overline R$.
 Se tiene por la continuidad de la funci\'{o}n corchete que $[x,y]= \lim [x_n,
 y_n]$. Seg\'{u}n la definici\'{o}n de rect\'{a}ngulo se sabe que $[x_n, y_n] \in
 R$, ya que $x_n,y_n \in R$. Entonces $[x,y]\in \overline R$ como
 se quer\'{\i}a.

 Sean ahora $x,y \in \mbox{int }R$. Entonces $[x,y]\in R$. Por la
 expansividad de $f$ se tiene que $x = [[x,y],x]$ e $y =
 [y,[x,y]]$.

 Sean $V_x, V_y$ entornos de $x$  e $y$ respectivamente, ambos
 contenidos en $R$. Por la continuidad de la funci\'{o}n corchete
 existe $V$ entorno de $[x,y]$ en $M$ tal que $[V,x] \subset V_x, \;\;
 [y,V] \subset V_y$. Entonces para todo $z \in V$ se cumple
 $$[[z,x],[y,z]] \in R$$ Siendo $\epsilon < \rho /4$ tenemos que $z= [[z,x],[y,z]]$
 y entonces $V \subset R$. Luego $[x,y]\in \mbox{int }R$ como se
 quer\'{\i}a.

  La intersecci\'{o}n no vac\'{\i}a de rect\'{a}ngulos es un rect\'{a}ngulo como se
  ve inmediatamente a partir de la definici\'{o}n de rect\'{a}ngulo. $\Box$

\vspace{.3cm}

{\bf \large 2.4 Definici\'{o}n de Partici\'{o}n de Markov }

\begin{definition}
\em Una  \em partici\'{o}n por cerrados \em de una variedad $M$ es un
cubrimiento finito de $M$ por cerrados $R_i$ con interiores dos a
dos disjuntos. El di\'{a}metro de la partici\'{o}n es el m\'{a}ximo de los
di\'{a}metros de los conjuntos cerrados que la componen.
\end{definition}
\begin{definition} \em \label{242}
Una partici\'{o}n  por cerrados de $M$ es una \em partici\'{o}n de Markov
\em para el difeomorfismo de Anosov $f$ si est\'{a} constituida por
rect\'{a}ngulos propios $R_1, R_2, \ldots, R_m$ y si para todo $ x \in
\mbox{int }R_i \cap f^{-1} \mbox{int }R_j$ se cumple:

\begin{itemize}
\item [i) ] $fW^s(x, R_i) \subset W^s(fx, R_j)$

\item [ii) ] $ fW^u(x, R_i) \supset   W^u(fx, R_j)$

\end{itemize}
\end{definition}

\begin{remark} \em \label{244}
La partici\'{o}n de Markov est\'{a} vinculada a $f$ a trav\'{e}s de la
definici\'{o}n de rect\'{a}ngulo y de las condiciones (i) y (ii). Si
$\mbox{int }R_i \cap f^{-1} \mbox{int }R_j \neq \emptyset$
entonces $f(R_i)$ se obtiene de $R_i$ al aplicarle $f$
comprimiendo las variedades $\epsilon$- estables y dilatando las
inestables.
\end{remark}

Sea una partici\'{o}n  ${\cal R}$ en $m$ subconjuntos de $M$(${\cal R
}$ no  es necesariamente de Markov),  con di\'{a}metro $\beta < \rho$
( $\rho $ es la constante de expansividad de $f$). Se cumple:
\begin{itemize}
\item[a) ] $\cap _{n \in \Z} f^{-n} R_{j_n} $ consta a lo sumo de
un punto (donde $j_n$ es una sucesi\'{o}n bi-infinita de n\'{u}meros en
$\{1,2, \ldots, m\}$). Esto es porque$$ x, y \in \bigcap _{n \in
\Z } f^{-n} R_{j_n} \;\; \Rightarrow \dist (f^n x, f^n y) <
\rho\;\; \forall n \in \Z \;\; \Rightarrow x= y$$
\item[b) ] Sea $\Sigma$ el conjunto de las sucesiones $\{j_n\}_{n \in \Z
}$ tales que $\cap _{n \in \Z} R_{j_n } \neq \emptyset$. Por lo
observado antes existe una funci\'{o}n $\Pi : \sigma \mapsto M$
definida por
$$\Pi (\{ j_n \}_{n \in \Z})= \bigcap _{n \in \Z} R_{j_n}$$
$\Pi $ es sobreyectiva pues toda \'{o}rbita $\{f^n x\}_{n \in \Z}$
est\'{a} cubierta por conjuntos de la partici\'{o}n. As\'{\i} dado $x \in M $
existe alguna sucesi\'{o}n $j_n$ tal que $f^n(x) \in R_{j_n } \; \;
\forall n \in \Z$. Luego $x \in \cap _{n \in \Z } f^{-n} R_{j_n}$.
\item[c) ] Si $x,y \in \cap _{n \geq 0 } f^{-n} R_{j_n}$ entonces
$\dist (f^n x, f^n y) \leq \beta \; \forall n \geq 0$. Luego $y
\in W^s _{\beta } (x) \cap R_{j_0 }$. Hemos probado que para
cualquier partici\'{o}n ${\cal R }$ se cumple:
$$  x \in \cap _{n \geq 0 } f^{-n} R_{j_n } \; \Rightarrow
\cap _{n \geq 0 } f^{-n} R_{j_n } \subset W^x _{\beta }(x) \cap
R_{j_0 }$$

\end{itemize}

\begin{remark}
\em Si $x = \Pi \{j_n\}_{n \in \Z }$ entonces $f^n (x) \in R_{j_n
} \forall n \in \Z$, o sea $f^n (fx) \in R_{n_{j+1}} \; \forall n
\in \Z$, de donde $fx = \Pi \{j_{n+1}\}_{n \in \Z }$.

Llamemos \em shift \em a la transformaci\'{o}n $\sigma : \Sigma
\mapsto \Sigma $ tal que a la sucesi\'{o}n $\{j_n\}_{n \in \Z }$ hace
corresponder la sucesi\'{o}n $\{j_{n+1}\}_{n \in \Z }$.

 Hemos obtenido que $\Pi \circ \sigma = f \circ \Pi$, es decir,
 conmuta el siguiente diagrama
 $$\begin{array}{ccccc}
    &  & \sigma &  &  \\
    & \Sigma & \rightarrow & \Sigma &  \\
   \Pi & \downarrow &  & \downarrow & \Pi \\
    & M & \rightarrow & M & \\
    & & f &  &  \\
 \end{array} $$
 Adem\'{a}s como $\Sigma$ y $f$ son invertibles se cumple que $\Pi \circ \Sigma ^n =
 f^n \circ \Pi \; \forall n \in \Z$. Luego:

 \em La funci\'{o}n sobreyectiva $\Pi$ lleva \'{o}rbitas del shift en
 \'{o}rbitas de $f$. \em
\end{remark}

Sea ${\cal R } = \{R_1, \ldots R_m \}$ una partici\'{o}n de Markov y
sea $x \in M$ un punto cuya \'{o}rbita por $f$ no corta a los bordes
de los cerrados $R_j$ de la partici\'{o}n, o sea $ x \in \cap _{n \in
\Z } f^{-n } (M \setminus \partial \cup _{j = 1 } ^{m} \partial
R_j)$. Sea $\{j_n\}_{n \in \Z }$ una sucesi\'{o}n tal que $f^n x \in
R_{j_n } \; \forall n \in \Z$ (o sea $\Pi (\{j_n\}_{n \in \Z })=
x$). Ahora, por construcci\'{o}n tenemos que $f^n x \in \mbox{int
}R_{j_n }$ y por definici\'{o}n de partici\'{o}n cerrada $\{j_n\}_{n \in
\Z }$ es \'{u}nica. La funci\'{o}n $\Pi $ es inyectiva sobre el conjunto
$\cap _{n \in \Z } f^{-n } (M \setminus  \cup _{j = 1 } ^{m}
\partial R_j)$. A continuaci\'{o}n veremos que ese conjunto es denso
en $M$ e invariante bajo $f$.

\vspace{.3cm}

{\bf \large 2.5 Borde de la partici\'{o}n de Markov }

\begin{definition}
\em Sea ${\cal R} = \{R_1, \ldots, R_m\}$ una partici\'{o}n de $M$ por
cerrados. Se llama \em Borde de ${\cal R}$ \em al conjunto
$$\partial {\cal R} = \bigcup _{j= 1}^m \partial  R_j$$ Si ${\cal R } $
es una partici\'{o}n de Markov se llama \em Borde estable de ${\cal
R}$ \em a $$\partial ^s {\cal R} = \cup _{j = 1 }^m \partial ^s
R_j$$ y se llama \em Borde inestable de ${\cal R}$ \em a
$$\partial ^u {\cal R} = \cup _{j = 1 }^m \partial ^u R_j$$ Se
observa de la proposici\'{o}n \ref{223} lo  siguiente:
$$\partial {\cal R } = \partial  ^s {\cal R} \cup \partial ^u {\cal R}$$
\end{definition}
\begin{proposition} \label{252}
\em Sea ${\cal R } = \{R_1, \ldots, R_m \}$ una partici\'{o}n por
cerrados de $M$. Entonces:
\begin{itemize}
\item [1) ] $\partial {\cal R}$ tiene interior vac\'{\i}o y es cerrado.
\item [2) ] $\cup _{j = 1} ^m \mbox {int } R_j$ es abierto y denso
en $M$.
\item [3) ] $A= M \setminus \cup _{n \in \Z } f^{-n } (\partial {\cal
R})$ es denso en $M$ e invariante por $f$.
\end{itemize}
\end{proposition}

{\em Demostraci\'{o}n: } (1) $\partial R_j$ es cerrado con interior
vac\'{\i}o. La uni\'{o}n finita de conjuntos en una variedad que son
cerrados con interior vac\'{\i}o es cerrada con interior vac\'{\i}o.

(2) Tomando el complemento $(\partial {\cal R}) ^c$ es abierto y
denso en $M$. Pero si $y \in (\partial {\cal R}) ^c$ entonces $y
\in \cap _{j = 1}^m (\partial {R}_j) ^c$. Como ${\cal R} $ cubre a
$M$ existe $j $ tal que $y \mbox{int } R_j$. Entonces $(\partial
{\cal R}) ^c \subset  \cup _{j = 1}^m \mbox {int } R_j$. Deducimos
que $\cup _{j = 1}^m \mbox {int } R_j$ es abierto y denso en $M$.

(3)$A = \cap _ {n \in \Z} (f^{-n} \partial {R}_j) ^c  = \cap _{n
\in \Z } f^{-n} ((\partial {R}_j) ^c)$ es denso en $M$ porque es
la intersecci\'{o}n numerable de abiertos densos. Adem\'{a}s $A$ es
invariante por $f$ porque $$f^{-1}(A ) = f^{-1} (\bigcap _{n \in
\Z } f^{-n} ((\partial {R}_j) ^c)) = \bigcap _{n \in \Z } f^{-n-1}
((\partial {R}_j) ^c) = A \;\; \Box$$

\vspace{.3cm}

{\bf \large 2.6 Propiedades de las particiones de Markov }

\vspace{.3cm}

La siguiente proposici\'{o}n permite aplicar las condiciones (i) y
(ii) de la definici\'{o}n \ref{242} de partici\'{o}n de Markov a otros
puntos $x \in M$ que no est\'{a}n necesariamente en $\mbox{int } R_i
\cap f^{-1}(\mbox {int }R_j)$.
\begin{proposition} \label{261}
Sea ${\cal R }= \{R_1, \ldots, R_m \}$ una partici\'{o}n de Markov de
$M$ para $f$. Si $ \mbox{int } R_i \cap f^{-1}(\mbox {int }R_j)
\neq \emptyset$ entonces para todo $y \in R_i \cap f^{-1}R_j $ se
cumple:
\begin{itemize}
\item [i) ] $f(W^s(y, R_i)) \subset W^s(fy, R_j)$
\item [ii) ] $f(W^u(y, R_i)) \supset W^u(fy, R_j)$

\end{itemize}
\end{proposition}

{\em Demostraci\'{o}n: } Sean $x \in \mbox{int } R_i \cap f^{-1}(\mbox
{int }R_j), \; \; y \in  R_i \cap f^{-1}(R_j)$. Es inmediato, a
partir de la definici\'{o}n de rect\'{a}ngulo y sabiendo que $x,y \in R_i$
lo siguiente:
$$W^s(y, R_i) = \{[y,z] : z \in W^s (x, R_i)\} = \{[y,z] : fz \in  f W^s (x, R_i)\}$$
Como $x \in \mbox{int } R_i \cap f^{-1}(\mbox {int }R_j)$ tenemos
por la condici\'{o}n (i) de la definici\'{o}n \ref{242} que se cumple:
$$f(W^s(x, R_i)) \subset W^s(fx, R_j) \subset R_j$$
Entonces $y,z \in R_i, \;\; fy, fz \in R_j$. As\'{\i} $\dist (y,z) \leq
\mbox{ diam } R_i \leq \delta < \epsilon $. An\'{a}logamente $\dist
(fy, fz) < \epsilon$, de donde $f[y,z] = [fy, fz ]$. Luego
$$fW^s (y, R_i) = \{ f[y,z]: z \in W^s (x, R_i)\}=
\{[fy,fz]: z \in W^s (x, R_i)\}=$$ $$= \{[fy,w]: w \in f W^s(x,
R_i)\} \subset \{ [fy, w]: w \in W^s (fx, R_j)\} = W^s(fy, R_j) $$
donde $fx, fy \in R_j$. Deducimos que  $fW^s (y, R_i) \subset W^s
(fy, R_j)$. An\'{a}logamente $fW^u (y, R_i) \supset W^u (fy, R_j)\;\;
\Box$

\begin{corollary}
Si $\{R_{j_n}\}_{n \geq 0 }$ es una sucesi\'{o}n de rect\'{a}ngulos de una
partici\'{o}n de Markov ${\cal R }$ con di\'{a}mtero $\beta  >0 $
suficientemente peque\~{n}o y tales que $\mbox { int } R_{j_n} \cap
f^{-1} \mbox { int }R_{j_{n+1}} \neq \emptyset \; \; \forall n
\geq 0$ entonces
$$x \in \bigcap _{n \geq 0} f^{-n} R_{j_n} \Rightarrow \bigcap _{n \geq 0 } f^{-n} R_{j_n}
= W^s _{\epsilon } (x, R_{j_0})$$
\end{corollary}

 {\em Demostraci\'{o}n: } Por lo observado en \ref{244} c) si se elige
 el di\'{a}metro de la partici\'{o}n de Markov menor que $\epsilon >0$
 obtenemos $\bigcap _{n \geq 0 } f^{-n} R_{j_n} \subset W^s_{\epsilon } (x,
 R_{j_0})$. Sea $y \in W^s_{\epsilon } (x, R_{j_0})$. Por
 \ref{261} se tiene $$fy \in W^s _{\epsilon }(fx, R_{j_1}), \;\;\;
 f^n y \in W^s_{\epsilon } (f^n y, R_{j_n}) \; \forall n \geq 0$$
Luego $y \in \bigcap _{n \geq 0 } f^{-n} R_{j_n}\; \Box$

\begin{proposition} Si ${\cal R}$ es una partici\'{o}n de Markov para
el difeomorfismo de Anosov $f$ entonces:
$$f(\partial ^s{\cal R}) \subset \partial ^s {\cal R}$$
$$f(\partial ^u{\cal R}) \supset \partial ^u {\cal R}$$
\end{proposition}

{\em Demostraci\'{o}n: } Sea $x \in \partial ^s {\cal R} = \bigcup
_{i= 1} ^m \partial ^s R_i$. Sea $i$ tal que $x \in \partial ^s
R_i$. En la proposici\'{o}n \ref{252} se prob\'{o} que $\bigcup _{j=1}^m
\mbox { int } R_j$ es denso en $M$. Entonces tambi\'{e}n lo es su
preimagen por el difeomorfismo $f$. Luego $$(f^{-1} \bigcup
_{j=1}^m \mbox{ int } R_j ) \bigcap \mbox { int }R_i
$$ es denso en $R_i$. Sea entonces $x_n \in (f^{-1} \bigcup
_{j=1}^m \mbox{ int } R_j ) \bigcap \mbox { int }R_i \; \forall n
$ tal que $x_n \rightarrow x$. Tenemos que $f(x_n) \in \bigcup
_{j=1}^m \mbox { int }R_j \; \forall n \geq 0 $, pero $j$ solo
puede tomar una cantidad finita de valores. Luego, existe una
subsucesi\'{o}n, que por comodidad seguimos llamando $x_n$, y un
\'{\i}ndice $j$ tal que $f(x_n) \in \mbox{ int }R_j  \;  \forall n \geq
0$.

El rect\'{a}ngulo $R_j$ es cerrado. Entonces $f(x) = \lim f (x_n) \in
R_j$. Tenemos entonces $$x_n \in (f^{-1}  \mbox{ int } R_j )
\bigcap \mbox { int }R_i$$ Luego por \ref{261} se cumple
$$f (W^u(x, R_i)) \supset W^u (fx, R_j)$$
Supongamos por absurdo que $fx \not \in \partial ^s R_j$. Existe
un entorno $V$ de $fx$ en $W^u_{\epsilon } (fx)$ contenido en $R_j
\bigcap W^u _{\epsilon } (fx)$. Entonces
$$f^{-1}(W^u_{\epsilon } (fx))\subset W^u_{\epsilon } (x)$$
As\'{\i} $x \in \mbox{int } W^u (x, R_i)$ en $W^u_{\epsilon }(x)$, o
sea, $x \not \in \partial R_i$ contra lo supuesto.

De igual forma se prueba que $f(\partial ^u {\cal R}) \supset
\partial ^u{\cal R}\; \; \Box$

\vspace{11cm}

\pagebreak

.

\pagebreak

\begin{center}
{\Large Particiones de Markov para difeomorfismos de Anosov - }
Eleonora Catsigeras
\end{center}

\section{Semiconjugaci\'{o}n con el  shift}

{\bf \large 3.1 Espacio de funciones bi-infinitas y funci\'{o}n shift
}

\vspace{.3cm}

Sea ${ P} = \{p_1, \ldots , p_m\}$ un conjunto finito de puntos en
$M$. Se denota con $P^{\Z}$ al espacio de las sucesiones
bi-infinitas de puntos en $P$. Tomando en $P$ la topolog\'{\i}a
discreta y en $P^{\Z}$ la topolog\'{\i}a producto asociada a ella, por
el teorema de Tychonov, el espacio $P^{\Z}$ es compacto y
metrizable.

Sea $q = \{q_j\}_{j \in \Z}, \; q_j \in P, \; q \in P^{\Z}$. Una
base local de abiertos en $q$ est\'{a} formada por los abiertos
$$I_N (q) = \{q' \in P^{\Z}: q_j = q'_j \,  \; \forall |j| \leq N \}$$
Una m\'{e}trica en $P^{\Z}$ est\'{a} dada por
$$\dist (q,q') = \sum _{n \in \Z} \frac {\dist (q_n, q'_n)}{2^{|n|}}$$

\begin{definition}
\em La funci\'{o}n o transformaci\'{o}n \em shift, \em denotada como
$\sigma : P^{\Z} \mapsto P^{\Z}$, es la trasnformaci\'{o}n  definida
por $\sigma (q) = q' $ donde $q' _n = q_{n+1}, \; \forall n \in
\Z$.
\end{definition}

Se observa que la funci\'{o}n shift $\sigma$ aplicada a $q$ consiste
en un corrimiento a la izquierda de los t\'{e}rminos de $q$: el mismo
t\'{e}rmino $q_0$ que antes ocupaba el lugar 0, despu\'{e}s de aplicarle
$\sigma$ ocupar\'{a} el lugar $-1$ (es decir es $q'_{-1}$),  el
t\'{e}rmino $q_1$ que antes ocupaba el lugar 1 pasar\'{a} a ocupar el
lugar 0 (es decir ser\'{a} $q' _0$) y as\'{\i} $q_j = q'_{j-1}$ para todo
$j \in \Z$. Es f\'{a}cil demostrar que $\sigma : P^{\Z} \mapsto
P^{\Z}$ es un homeomorfismo.

 \vspace{.3cm}

{\bf \large 3.2 Semiconjugaci\'{o}n }

\vspace{.3cm}

Sean $M, M'$ dos espacios topol\'{o}gicos, y sean $f, f'$ dos
homeomorfismos en $M$ y $M'$ respectivamente.
\begin{definition}
\em Una funci\'{o}n $\theta: M' \mapsto M$ se llama \em
semiconjugaci\'{o}n \em de $f$ con $f'$ si cumple:
\begin{itemize}
\item [i) ] $\theta$ es continua y sobreyectiva
\item [ii) ] $f \circ \theta = \theta \circ f'$
\end{itemize}
\end{definition}
Se observa que
$$f \circ \theta = \theta \circ f'\; \; \; \Rightarrow \; \; \;
f  ^n \circ \theta = \theta \circ f'^n \; \; \forall n \in \Z$$
Luego, toda \'{o}rbita en $M'$ seg\'{u}n $f'$ es llevada por $\theta $ a
alguna \'{u}nica \'{o}rbita por $f$ en $M$ y toda \'{o}rbita en $M$ por $f$ es
corresponde a alguna (no necesariamente \'{u}nica) \'{o}rbita por $f'$ en
$M'$.
\begin{definition}
\em Una semiconjugaci\'{o}n se llama \em conjugaci\'{o}n \em entre $f$ y
$f'$ si es un homeomorfismo.
\end{definition}
{\bf \large 3.3 Semiconjugaci\'{o}n de los difeomorfismos de Anosov
con el shift}

\vspace{.3cm}

Sea $\beta >0 $ arbitrario dado. Sea $f: M \mapsto M $ un
difeomorfismo de Anosov. Por el teorema \ref{135} el difeomorfismo
$f$ es topol\'{o}gicamente estable. Elijamos $\alpha >0 $ tal que toda
pseudo-\'{o}rbita de $f$ est\'{a} $\beta $ acompa\~{n}ada por una \'{o}rbita de
$f$.

Sea $0< \gamma < \min (\beta, \alpha /2)$ tal que $$\dist (x,y )<
\gamma \; \Rightarrow \; \dist (fx, fy) < \alpha /2$$ Tal n\'{u}mero
$\gamma $ existe porque $f$ es continua en $M$ compacta.

Siendo $M$ compacta existe un conjunto finito $P = \{p_1, \ldots,
p_m\}$ de puntos de $M$, centros de bolas de radio $\gamma$ que
cubren $M$. Dado $x \in M$ existe $p_j \in P$ tal que $\dist (x,
p_j) < \gamma$. Es decir $P$ es un conjunto $\gamma$-denso en $M$.

Sea $\Sigma (P) = \{q \in P^{\Z}: \dist (fq_j, f q _{j+1})< \alpha
\; \forall j \in \Z \}$

$\Sigma (P)$ es el conjunto de las $\alpha $-pseudo-\'{o}rbitas de $f$
que est\'{a}n formadas con puntos de $P$.

Si adem\'{a}s elegimos $\beta < \rho /2$, donde $\rho$ es la constante
de expansividad de $f$, se cumple, en virtud de la estabilidad
topol\'{o}gica de $f$ dada por el teorema \ref{135}, lo siguiente:

Para todo $q \in \Sigma (P)$ existe un \'{u}nico $\theta (q) \in M$
tal que $$dist (f^n(\theta (q)), q_n) \leq \beta \; \; \forall n
\in Z
$$
\begin{lemma}
La aplicaci\'{o}n $\theta : \Sigma (P) \mapsto M$  es sobreyectiva.
\end{lemma}

{\em Demostraci\'{o}n: } Sea $x \in M$. Demostremos que existe alg\'{u}n
$q \in \Sigma (P)$ tal que $x = \theta (q)$.

La \'{o}rbita $\{f^n(x)\}_{n \in \Z}$ se puede aproximar por $q =
\{q_n\}_ {n \in \Z} \in P^{\Z}$ de modo que $$\dist (f^n(x), q_n)
\leq \gamma \; \; \forall n \in \Z$$ porque $P$ es $\gamma$-denso
en $M$.

Entonces $\dist (fq_n, q_{n+1}) \leq \dist (fq_n, f^{n+1} x) +
\gamma$. De acuerdo a la elecci\'{o}n de $\gamma$, siendo $\dist (f^n
x, q_n) < \gamma$, se cumple $\dist (f^{n+1} x, f q_n)< \alpha
/2$.

As\'{\i} $\dist (f q_n, q_{n+1}) \leq \alpha /2 + \gamma$.

Siendo $\gamma < \alpha /2$ se cumple $\dist (f q_n, q_{n+1}) <
\alpha$, o sea $\{q_n\}$ es una $\alpha$-pseudo-\'{o}rbita de $f$.
Luego $q \in \Sigma (P)$.

Como $\{f^nx\}_{n \in \Z}$ es una \'{o}rbita que $\gamma $-acompa\~{n}a a
$q$ por construcci\'{o}n, y siendo $\alpha < \beta$ resulta $\dist
(f^n x, q_n) < \beta \; \forall n \in Z$. Entonces $x = \theta
(q)$ como se quer\'{\i}a demostrar. $\Box$

\begin{lemma}
La aplicaci\'{o}n $\theta : \Sigma (P) \mapsto M$ es continua
\end{lemma}

{\em Demostraci\'{o}n: } Sea $q^n \rightarrow q \in \Sigma (P), \; \;
\theta (q^n) = x_n \in M$. La sucesi\'{o}n $x_n$ puede suponerse
convergente $x_n \rightarrow x_0 $debido a la compacidad de la
variedad $M$.

Se tiene que $\dist (q^n_j, f^j (x_n)) \leq \beta \; \; \forall j
\in \Z$ por la construcci\'{o}n de la funci\'{o}n $\theta$.

Sea $j \in \Z $ fijo. Como $q^n \rightarrow q \in \Sigma (P)$,
existe $N(j)$ tal que para todo $n > N(j)$ se cumple $q^n \in I_j
(q)$, es decir $q^n _i = q_i \; \forall |i| \leq j$.

De lo anterior se deduce que para todo $n > N(j)$:
$$ \dist (q_j, f^j (x_n)) \leq \beta$$
Tomando $n \rightarrow \infty$, en virtud de la continuidad de $f$
se deduce que
$$ \dist (q_j, f^j (x_0)) \leq \beta \; \; \forall j \in \Z$$
Entonces por construcci\'{o}n de la funci\'{o}n $\theta$ se cumple que
$x_0 = \theta (q)$ y luego $\theta (q^n ) \rightarrow \theta
(q).\; \; \Box$

\begin{theorem} \label{334}
La funci\'{o}n $\theta : \Sigma (P) \mapsto M$ es una semiconjugaci\'{o}n
del difeomorfismo de Anosov $f: M \mapsto M$ con el shift $\sigma
$  restringido a $\Sigma (P)$
\end{theorem}

{\em Demostraci\'{o}n: } La funci\'{o}n shift $\sigma $, cuando
restringida a $\Sigma (P)$,  tiene codominio en $\Sigma (P)$ pues
$\Sigma (P)$ es invariante bajo $\sigma$. En otras palabras $\left
. \sigma \right |
 _{\Sigma (P)}: \Sigma (P) \mapsto \Sigma (P)$.

 Para demostrar que $\theta$ es una semiconjugaci\'{o}n entre $f$ y $\left
. \sigma \right |
 _{\Sigma (P)}$ alcanza demostrar que el diagrama siguiente
 conmuta, pues ya se sabe que $\theta$ es continua y sobreyectiva:

 $$\begin{array}{ccccc}
     &   & \sigma  &   &   \\
     & \Sigma (P) & \mapsto & \Sigma (P) &  \\
   \theta & \downarrow &   & \downarrow & \theta \\
     & M &  \mapsto & M &   \\
    &  & f &  &  \\
 \end{array} $$
Sea $q \in \Sigma (P)$. Sean $x = \theta (q), \; q' = \sigma (q)$.
Alcanza probar que $f(x) = \theta (q')$. Sea $x' = \theta (q')$.
Entonces  $$\dist (f^n(x'), q'_n) \leq \beta \; \; \forall n \in
\Z$$ Pero $q'_n = q_{n+1}$ pues $q' = \sigma (q)$. Entonces
$$\dist (f^n(x'), q_{n+1}) \leq \beta \; \; \forall n \in
\Z$$
$$\dist (f^{n+1}(f^{-1} (x')), q_{n+1}) \leq \beta \; \; \forall n \in
\Z$$ Luego $f^{-1}(x') = \theta (q) = x$, de donde $ x' = f(x)$
como quer\'{\i}amos probar. $\; \Box$

\vspace{.5cm}

{\bf \large 3.4 Conjuntos estable e inestable en el espacio de
sucesiones}

\vspace{.3cm}

Sea $q \in \Sigma (P)$ donde $\Sigma (P)$ es el subconjunto de
$P^{\Z}$ (sucesiones bi-infinitas) definido en la secci\'{o}n 3.3

\begin{definition}
\em Se llama \em conjunto estable \em por $q$ en $\Sigma (P)$ a:
$$\widehat W^s(q) = \{q' \in \Sigma (P): \dist (\sigma q, \sigma q') \rightarrow
 _{n \rightarrow + \infty} 0\}$$

\em Se llama \em conjunto inestable \em por $q$ en $\Sigma (P)$ a:
$$\widehat W^u(q) = \{q' \in \Sigma (P): \dist (\sigma q, \sigma q') \rightarrow
 _{n \rightarrow - \infty} 0\}$$
\end{definition}
\begin{remark}
\em Se sabe que $dist (q, q') = \sum _{j \in \Z} \dist (q_j, q'_j)
/2^{|j|}$. Luego: $$dist (\sigma ^n q, \sigma ^n q') = \sum _{j
\in \Z} \dist (q_{n+j}, q'_{n+j}) /2^{|j|}$$ Si $\dist (\sigma ^n
q, \sigma ^n q') \rightarrow 0$ entonces existe $N$ tal que
$\forall n >N$ se cumple
$$\sum _{j \in \Z} \frac {\dist(q_{n+j}, q' _{n+j})}{2^{|j|}} < \min \{\dist (p_i,
p_j): i \neq j, p_i, p_j \in P\}$$ Lo anterior se cumple si y solo
si $q_n = q' _n$ para todo $n$ suficientemente grande. Luego:
$$\widehat W^s (q) = \{q' \in \Sigma (P): q'_n= q_n\; \forall n \mbox{ suficientemente grande}
\}= $$  $$\widehat W^s (q) = \bigcup_{N \geq 0} \{q' \in \Sigma
(P): q'_n= q_n\; \forall n \geq N \}$$ An\'{a}logamente:
$$\widehat W^u (q) = \bigcup_{N \leq 0} \{q' \in \Sigma (P):
q'_n= q_n\; \forall n \leq N \} $$
\end{remark}
\begin{definition}
\em Se llama \em $0-$conjunto estable (e inestable) \em por $q$ a
$$\widehat W^s_0(q)= \{q' \in \Sigma (P): q'_n = q_n \; \forall n \geq 0\}$$
(respectivamente a:
$$\widehat W^u_0(q)= \{q' \in \Sigma (P): q'_n = q_n \; \forall n \leq 0\} )$$

\end{definition}

\begin{remark}
\em
\begin{itemize}
\item[1) ] $\widehat W^s_0(q) \subset \widehat W^s(q), \; \; \;
\widehat W^u_0(q) \subset \widehat W^u(q)$
\item[2) ] $\sigma \widehat W^s_0(q) \subset \widehat W_0 ^s(\sigma q), \; \; \;
 \sigma \widehat W^u_0(q) \supset \widehat W^u_0 (\sigma q)$
\item[3) ] Si $q, q'\in \Sigma (P)$ y si $q_0 = q'_0$, entonces
puede construirse una \'{u}nica $q"$ tal que $q"_j = q_j \; \forall j
\geq 0, \; \; \; q"_j = q'_j \; \forall j \leq 0$. Se cumple
$$q" = \widehat W^s_0 (q) \bigcap \widehat W^u _0 (q')$$

\end{itemize}
\end{remark}
\begin{definition}
La funci\'{o}n corchete en el espacio de sucesiones es:
$$[,]: \{(q,q') \in (\Sigma (P))^2: q_0 = q' _0\} \mapsto \Sigma (P)$$
definida por $$[q,q'] = \widehat W _0 ^s (q) \bigcap \widehat W
^u_0 (q')$$ o sea $q" = [q,q'] $ si y solo si $q"_j = q_j  \;
\forall j \geq 0, \; \; \; q"_j = q'_j \; \forall j \leq 0$.

\end{definition}

\begin{proposition} \label{346}
\em Sea $\theta $ la semiconjugaci\'{o}n definida en la secci\'{o}n 3.3.
Si $\beta $ es suficientemente peque\~{n}o entonces
$$\theta [q,q'] = [\theta q, \theta q']\; \forall q, q' \in \Sigma (P) \mbox { tales que }
q_0 = q' _0$$
\end{proposition}
{\em Demostraci\'{o}n: } Sean $q, q' \in \Sigma (P)$ tales que $q_0 =
q' _0$. Llamemos $x = \theta (q), \; \; y = \theta (q'), \; \; z =
\theta [q,q']$. Hay que demostrar que $z = [x,y]$ o sea que $z =
W^s_{\epsilon }(x) \cap W^u_{\epsilon }(y)$.

Por construcci\'{o}n y por definici\'{o}n de la funci\'{o}n $\theta $ se
cumple:
$$\dist (f^n x, q_n) \leq \beta \; \; \; \dist (f^n y, q'_n) \leq \beta \; \; \; \forall n \in \Z
$$ $$\dist (f^n z, q_n) \leq \beta \; \forall n \geq 0, \; \; \dist (f^n z, q' _n ) \leq
\beta \; \forall n \leq 0 $$ porque $z = \theta [q,q']$.

Entonces $$\dist (f^n x, f^n z) \leq 2 \beta \; \forall n \geq 0,
\; \; \; \dist (f^n y, f^n z) \leq 2 \beta \; \forall n \leq 0$$
Eligiendo $2 \beta < \min (\epsilon, \delta /2)$ se tiene que $z
\in W^s_{\epsilon } (x) \cap W^u_ {\epsilon } (y)$. Adem\'{a}s $\dist
(x,y) \leq \dist (x,z) + \dist (y,z) \leq 4 \beta < \delta$
Entonces por la expansividad de $f$, el punto $z \in M$ es el
\'{u}nico en $W^s_{\epsilon } (x) \cap W^u_ {\epsilon } (y). \; \; \;
\; \Box$

La semiconjugaci\'{o}n definida en la secci\'{o}n 3.3 conmuta con la
funci\'{o}n corchete.

\vspace{.5cm}

{\bf \large 3.5 Construcci\'{o}n de un cubrimiento con rect\'{a}ngulos}

\begin{proposition} \label{351}
Sea $\beta >0 $ dado suficientemente peque\~{n}o, como en la
Proposici\'{o}n \ref{346}. Si $p_s \in P$ (seg\'{u}n la secci\'{o}n 3.3)
entonces $T_s = \theta \{q \in \Sigma (P): q_0 = p_s\} $ es un
rect\'{a}ngulo cerrado de $M$ con di\'{a}metro a lo sumo $2 \beta $.
\end{proposition}

{\em Demostraci\'{o}n: } Si $q, q'$ cumplen $q_0 = q' _0 = p _s$
entonces tambi\'{e}n se cumple por construcci\'{o}n de $[q,q']$ que
$[q,q']_0 = p _s$.

Para demostrar que $T_s$ es un rect\'{a}ngulo en $M$ alcanza tomar
$x,y \in T_s$ y demostrar que $[x,y] \in T_s$. Pero $x,y \in T_s
 \Rightarrow \theta (q) = x, \theta (q') = y $ con $q_0 = q' _0 =
 p_s$.
 Luego seg\'{u}n \ref{346} se tiene $[x,y] = \theta [q,q']$. Entonces
 $[x,y]\in T_s$ (por construcci\'{o}n de $T_s$).

 Adem\'{a}s $\dist (x,y) < \dist (x, q_0) + \dist (q'0, y) < 2 \beta$.
 Entonces $\mbox{diam} T_s \leq 2 \beta$.

 $T_s$ es cerrado porque es la imagen continua de $\{q \in \Sigma (P): q_0 = p_s\}$
 que es compacto en $P^{\Z}$ (ya que es cerrado en el espacio compacto
 $P^{\Z}$).  $\Box$

 \begin{corollary}
 La familia de rect\'{a}ngulos $\tau = \{T_i\}_{i = 1, \ldots, m}$ (construidos
 seg\'{u}n la proposici\'{o}n
 \ref{351}) es un cubrimiento finito de $M$ por rect\'{a}ngulos
 cerrados de di\'{a}metro a lo sumo $2 \beta$.
 \end{corollary}

 {\em Demostraci\'{o}n: } Alcanza ver que $\cup _{i=1} ^m T_i = M$.
 Como $\theta$ es sobreyectiva, todo punto $x \in M$ es $x = \theta
 (q)$ con $q \in \Sigma (P)$. Sea $q_0 \in P = \{p_1, \ldots,
 p_m\}$. Entonces existe un sub\'{\i}ndice $s = 1, 2, \ldots, m$ tal
 que $q_0 = p_s$, o sea $x \in T_s$. Luego $M = \cup _{i=1} ^m T_i $
 como se quer\'{\i}a probar. $\Box$

\vspace{.5cm}

{\bf \large 3.6 Propiedades del cubrimiento por rect\'{a}ngulos}

\begin{proposition} \label{361}
Si $x = \theta (q)$ con $q_0 = p_s$ entonces $\theta (\widehat W_0
^s q) = W^x (x, T_s)$
\end{proposition}

{\em Demostraci\'{o}n: } Si $y \in \theta (\widehat W_0 ^s q) $,
entonces $y = \theta (q')$ con $q'_j = q_j , \forall j \geq 0$.
As\'{\i} $[\theta q, \theta q'] = [x,y] = \theta [q,q']$ (por la
Proposici\'{o}n \ref{346}).

Como $q_j' = q_j \; \forall j \geq 0$ tenemos que $[q,q']= q'$.

Entonces $[x,y] = \theta q' = y$, de donde $y \in W^s(x, T_s)$.

Hemos probado que $\theta (\widehat W_0 ^s q) \subset W^x (x,
T_s)$. Rec\'{\i}procamente, si $y \in W^s(x, T_s) \subset T_s$ entonces
existe $q'$ con $q'_0 = p_s$ tal que $y = \theta (q')$ (por
construcci\'{o}n del rect\'{a}ngulo $T_s$).

Adem\'{a}s si $y \in W^s(x, T_s)$ entonces $y = [x,y]$.

Aplicando la proposici\'{o}n \ref{346} se obtiene:

$$y = [x,y] = [\theta q, \theta q'] = \theta [q,q'] \in \theta (\widehat W^s_0 q)$$
Hemos probado  entonces que $W^s (x, T_s) \subset \theta (\widehat
W^s_0 q). \; \; \Box$

\begin{remark}
\em La proposici\'{o}n anterior caracteriza las $\epsilon$-variedades
estables en los rect\'{a}ngulos $T_s$ de $M$: la semiconjugaci\'{o}n
$\theta$ lleva $0-$ variedades estables en el espacio $\Sigma (P)$
en $\epsilon$- variedades estables en rect\'{a}ngulos $T_s$ de $M$.
\end{remark}

La siguiente proposici\'{o}n ser\'{a} utilizada en la secci\'{o}n 5 de este
trabajo para demostrar el teorema de Sinai (existencia de una
partici\'{o}n de Markov para $f$, difeomorfismo de Anosov).

\begin{proposition}
Sea $q \in \Sigma (P)$ con $q_0 = p_s, q_1 = p_t$ (seg\'{u}n la
definici\'{o}n al principio de la secci\'{o}n 3.3).  Sea $x = \theta (q)$.
Entonces:

\begin{itemize}
\item [ i) ] $f W^s(x, T_s) \subset W^s(fx, T_t)$
\item [ii) ] $f W^u(x, T_s) \supset W^u(fx, T_t)$
\end{itemize}
\end{proposition}

{\em Demostraci\'{o}n: } Se tiene $fx = \theta (\sigma q), \; (\sigma
q)_0 = q_1 = p _t , \; \; x = \theta (q), \; q_0 = p_s$.

Por la Proposici\'{o}n \ref{361} tenemos que $W^s(x, T_s) = \theta
(\widehat W^s_0 q), \; \; W^s(fx, T_t)= \theta (\widehat W^s_0
(\sigma q)$.

Por el Teorema \ref{334} $f W^s(x, T_s) = f \circ \theta (\widehat
W ^s _0 q) = \theta \sigma (\widehat W^s_0 q)$.

De la definici\'{o}n de $\widehat W ^s _0 (q)$  y de la definici\'{o}n de
la funci\'{o}n shift $\sigma$, es inmediato que $\sigma (\widehat
W^s_0 q )\subset \widehat W^s _0 (\sigma q))$. Entonces:
$$f W^s(x, T_s) \subset \theta (\widehat W^s_0 \sigma q) = W^s(fx, T_t))$$
En forma similar, utilizando conjuntos inestables en vez de
estables, se prueba (ii). $\; \Box$

\begin{remark}
\em Este procedimiento ha permitido construir un cubrimiento $\tau
= \{T_1, \ldots, T_m\}$ de $M$ por rect\'{a}ngulos cerrados de
di\'{a}metro menor que un n\'{u}mero positivo dado y que cumplen las
condiciones (i) y (ii) de la proposici\'{o}n anterior. Estas
condiciones son similares a las exigidas en la definici\'{o}n de
partici\'{o}n de Markov en el par\'{a}grafo \ref{242}.

El cubrimiento $\tau$ no es necesariamente una partici\'{o}n de Markov
porque los interiores de los rect\'{a}ngulos de $\tau$ no son en
general disjuntos dos y a dos y los rect\'{a}ngulos no son
necesariamente propios. A partir del cubrimiento $\tau$, que
cumple (i) y (ii), refin\'{a}ndolo apropiadamente, se construir\'{a} una
partici\'{o}n de Markov.
\end{remark}

\vspace{3cm}

\pagebreak

\begin{center}
{\Large Particiones de Markov para difeomorfismos de Anosov - }
Eleonora Catsigeras
\end{center}

\section{M\'{e}todo constructivo de la partici\'{o}n}

En la secci\'{o}n 3.5 se construy\'{o} un cubrimiento finito $\tau $ de la
variedad $M$, con rect\'{a}ngulos cerrados $\{T_1, T_2, \ldots, T_m\}$
para el difeomorfismo de Anosov $f$, que cumplen las condiciones
(i) y (ii) de la Definici\'{o}n \ref{242} de Partici\'{o}n de Markov.

En esta secci\'{o}n se refinar\'{a} el curbrimiento $\tau$ para obtener
ahora una  partici\'{o}n ${\cal R}$ por cerrados de $M$ (con
interiores dos a dos disjuntos) que sean adem\'{a}s conjuntos propios
(cada cerrado es la adherencia de su interior).

Finalmente en la secci\'{o}n 5 se demostrar\'{a} que esa partici\'{o}n ${\cal
R}$ es una partici\'{o}n de Markov para $f$.

\vspace{.5cm}

{\bf \large 4.1 Primer refinamiento del cubrimiento}

\vspace{.3cm}

A partir del cubrimiento $\{T_i\}_{i= 1, 2, \ldots, m} = \tau$
definamos otro cubrimiento m\'{a}s fino, de la siguiente forma:

\begin{definition}
\em Sean $T_j, T_k \in \tau$. Se definen:
\begin{itemize}
\item $T_{jk}^1 = \{x \in T_j: W^u(x, T_j) \cap T_k \neq \emptyset, \;
W^s(x, T_j) \cap T_k \neq \emptyset\}$
\item $T_{jk}^2 = \{x \in T_j: W^u(x, T_j) \cap T_k \neq \emptyset, \;
W^s(x, T_j) \cap T_k = \emptyset\}$
\item $T_{jk}^3 = \{x \in T_j: W^u(x, T_j) \cap T_k = \emptyset, \;
W^s(x, T_j) \cap T_k \neq \emptyset\}$
\item $T_{jk}^4 = \{x \in T_j: W^u(x, T_j) \cap T_k = \emptyset, \;
W^s(x, T_j) \cap T_k = \emptyset\}$
\end{itemize}
\end{definition}
\begin{remark}
\em
\begin{itemize}
\item[(1)] $T_j $ es la uni\'{o}n disjunta $\cup _{n=1}^4 T_{jk}^n$.
\item[(2)] Si $n_1 \neq n_2$ entonces $\mbox{int}T_{jk}^{n_1} \cap T_{jk}^{n_2} =
\emptyset = T_{jk}^{n_1} \cap T_{jk}^{n_2}$.
\item[(3)] $T_{jk}^1 = T_j \cap T_k$. En efecto: $T_j \cap T_k \subset
T^1_{jk}$. Adem\'{a}s si $x \in T_{jk}^1$ entonces existen $y,z \in
T_j \cap T_k$ tales que $y \in W^u_{\epsilon} (x), \; z \in
W^s_{\epsilon }(x)$. Luego $x = [ x,y]$. Pero por definici\'{o}n de
rect\'{a}ngulo, como $y,z \in T_j \cap T_k$, entonces $x = [z,y] \in
T_j \cap T_k$.
\end{itemize}
\end{remark}

\begin{proposition}
\em Si $T_{jk}^n \neq \emptyset$ entonces $T_{jk}^n$ es un
rect\'{a}ngulo.
\end{proposition}

{\em Demostraci\'{o}n: } Sean $x,y \in T_{jk}^n \subset T_j$. Entonces
$z \in [x,y] \in T_j$ porque $T_j$ es un rect\'{a}ngulo.

Por la proposici\'{o}n \ref{215}, como $z \in W^s(x, T_j)$, tenemos
que $W^s (z, T_j) = W^s (x, T_j)$. An\'{a}logamente $W^u(z, T_j) = W^u
(y, T_j)$.  Entonces $W^s (z, T_j)$ y $W^u (z, T_j)$ cortan a $T_k
$ si y solo si lo hacen $W^s(x, T_j)$ o respectivamente $W^u(y,
T_j)$.  Luego $z \in T^n_{jk}$ como quer\'{\i}amos. $\; \Box$

\vspace{.5cm}

{\bf \large 4.2 Segundo refinamiento de la partici\'{o}n}

\vspace{.3cm}

Los rect\'{a}ngulos $T_{jk}^n$ construidos al principio de la secci\'{o}n
4.1 cubren a $M$ pero no tienen necesariamente interiores
disjuntos. Tampoco son todos propios porque no son todos cerrados.
Construiremos un refinamiento ${\cal R}$ de $\{T_{jk }^n\}$.

En primer lugar hay que observar que un punto $x \in M$ puede
pertenecer a varios rect\'{a}ngulos de la familia $\{T_{jk }^n\}$.

\begin{definition}
\em Dado $x \in M$ sea
$$H(x) = \{(j,k,n): x \in T_{j,k}^n\}$$
$$R(x) = \bigcap _{(j,k,n) \in H(x) } \mbox{ int } \overline T_{j,k}^n$$
donde $R(x)$ podr\'{\i}a ser vac\'{\i}o.
$$Z^* = \{x \in M: x \in \mbox { int } T_{j,k}^n \; \; \forall (j,k,n) \in H(x))\}$$
$${\cal R} = \{ \overline {R(x)}\, \}_{x \in Z^*}$$
\end{definition}

\begin{remark}
\em
\begin{itemize}
\item[(1)] ${\cal R}$ es una familia finita, porque $H(x)$ es un
subconjunto del conjunto finito de todos los posibles \'{\i}ndices
$\{(j,k,n): 1 \leq j \leq m, 1 \ leq k \leq m, \ leq n \leq 4\}$
\item[(2)] $ x \in Z^*$ si y solo si toda vez que $x \in
T_{j,k}^n$ se cumple $x \in \mbox{ int } T_{j,k}^n$.
\item[(3)] Si $x \in Z^*$ entonces $x \in R (x)$; luego en ese
caso $R(x) \neq \emptyset$.
\item[(4)] Si $x \in Z^*$ entonces $R(x)$ es un rect\'{a}ngulo
abierto, porque no es vac\'{\i}o y es intersecci\'{o}n finita de
rect\'{a}ngulos abiertos.
\item[(5)] Si $x \in Z^* $ entonces $\overline {R(x)}$ es un
rect\'{a}ngulo propio pues
$$\overline{R(x)} \supset \overline{\mbox{ int } \overline{R(x)}}
\supset \overline{\mbox{ int } {R(x)}} =  \overline{ {R(x)}}
$$ pues $R(x)$ es abierto. Entonces $\overline
{R(x)} =  \overline{\mbox{ int } \overline{R(x)}}$.
\end{itemize}
\end{remark}

De las observaciones anteriores se deduce que ${\cal R}$ es una
familia finita de rect\'{a}ngulos propios que cubren $Z^*$. Probaremos
que ${\cal R}$ es una partici\'{o}n de $M$ (m\'{a}s a\'{u}n ser\'{a} una partici\'{o}n
de Markov), para lo cual demostraremos que:

\begin{itemize}
\item [ I) ] La familia ${\cal R}$ cubre a $M$ lo que se
demostrar\'{a} en la secci\'{o}n 4.4.
\item [II) ] Dos rect\'{a}ngulos de la familiar ${\cal R}$ que sean
diferentes tienen interiores disjuntos (esto se demostrar\'{a} en la
secci\'{o}n 4.4).
\end{itemize}

Las afirmaciones (I) y (II) se ´probar\'{a}n a partir de las
definiciones de los rect\'{a}ngulos $T_{j,k}^n$ y del conjunto $Z^*$.
En especial la densidad de $Z^*$ en $M$ que se demuestra a
continuaci\'{o}n juega un papel importante en la prueba.

\vspace{.5cm}

{\bf \large 4.3 Densidad del conjunto $Z^*$ cubierto por la
partici\'{o}n.}

\begin{proposition}
El conjunto $Z^*$ definido en la secci\'{o}n \em 4.2 \em es abierto y
denso en $M$. Adem\'{a}s \em
$$Z^* = \bigcap _{j,k = 1}^m  \left
( \bigcup _{n= 1}^4 \mbox{ int }T_{j, k}^n \cup T_j^c \right )$$
\end{proposition}

{\em Demostraci\'{o}n: }

Por definici\'{o}n
$$Z^* = \{x \in M: x \in \mbox{ int } T_{jk}^n \; \forall (j,k,n) \in H(x)\}
 = \{ x \in M: x \in T_{jk}^n \Rightarrow \mbox{ int } T_{jk}^n\}$$
 Se prob\'{o} que $M = \cap _{j=1}^m T_j$. Adem\'{a}s por construcci\'{o}n $T_j
 $ es la uni\'{o}n disjunta $ \cup _{n=1}^4 T_{j,k}^n \; \forall
 (j,k)$.

 Fijados $j,k (1 \leq j,k \leq m)$ y dado un punto $x \in M$ se
 cumple $x \in T_j^c \cup \cup_{n=1}^4 T_{jk}^n$. Si el punto $x \in
 Z^*$ entonces $ x \in T_j^c \cup \cup_{n=1}^4 T_{jk}^n \; \forall j, k$

Y rec\'{\i}procamente, si $x \in T_j^c \cup \cup_{n=1}^4 T_{jk}^n \;
\forall j, k$ entonces
$$x \in \{x \in M: x \in T_{jk}^n \Rightarrow x \in \mbox{ int } T_{jk}^n\} = Z^*$$

Luego $$Z^* = \bigcap _{j,k = 1}^m \left ( T_j^c \bigcup
\bigcup_{n=1}^4 T_{jk}^n \right)$$

Entonces $Z^*$ es abierto, por se intersecci\'{o}n finita de abiertos.

Para demostrar que $Z^* $ es denso en $M$ alcanza probar que cada
uno de los siguientes abiertos
$$Z_{jk}^* = T_j^c \bigcup \bigcup _{n=1}^4 T_{jk}^n$$
es denso en $M$.

Para eso es suficiente tomar un abierto cualquiera no vac\'{\i}o
contenido en $T_j$ y probar que corta a $\bigcup _{n=1}^4 \mbox{
int } T_{jk}^n$.

El lema que sigue permite caracterizar $\bigcup_{n=1}^4\mbox{ int
} T_{jk}^n$ como intersecci\'{o}n de dos abiertos densos en $\mbox {
int } T_j$. Luego por el teorema de Baire, $\bigcup_{n=1}^4\mbox{
int } T_{jk}^n$ ser\'{a} denso en $\mbox { int } T_j$ como queremos.

\begin{lemma}
Dados $j,k$ sean \em
$$A _{j,k} = \mbox{ int } \{x \in T_j: W^s(x,T_j) \cap T_k = \emptyset\} \; \cup
 \;  \mbox{ int } \{x \in T_j : W^s (x, T_j) \cap T_k \neq \emptyset\}$$
 $$B_{j,k} =   \mbox{ int } \{x \in T_j: W^u(x,T_j) \cap T_k = \emptyset\} \;
 \cup \;
 \mbox{ int } \{x \in T_j : W^u (x, T_j) \cap T_k \neq \emptyset\}$$
\em Entonces
\begin{itemize} \em
\item[ 1) ] $A_{j,k} \cap B_{j,k}  = \cup_{n=1}^4 \mbox { int } T_{j,k}^n$
\item[ 2) ] $A_{j,k}$ y $B_{j,k}$ \em son abiertos y densos en \em $\mbox{ int }T_{j}$
\end{itemize}
\end{lemma}
{\em Demostraci\'{o}n: } 1) Teniendo en cuenta las definiciones en la
secci\'{o}n 4.1 se tiene que
$$A_{j,k} \cap B_{j,k} \cap T{j,k}^n = \mbox {int } T_{j,k}^n \; \forall j,k,n$$
Entonces
$$A_{j,k} \cap B_{j,k} \cap  \cup _{n=1}^4 T_{j,k}^n = \cup _{n=1}^4 \mbox{ int } T_{j,k} ^n
\; \; \forall j,k $$ Sabemos que $\cup_{n=1}^4 T_{j,k}^n = T_j \ni
A_{j,k} \cap B_{j,k} $ y entonces $A_{j,k} \cap B_{j,k} = \cup
_{n=1}^4 \mbox { int } T_{j,k}^n \; \forall j,k$

2) $A_{j,k}$ y $B_{j,k}$ son abiertos por construcci\'{o}n. Mostremos
que $A_{j,k}$ es denso en $\mbox{ int } T_j$. Sea $V$ un abierto
no vac\'{\i}o contenido en $T_j$. Basta probar que $V$ contiene alg\'{u}n
punto de $A_{j,k}$.
$$V = \{ y \in V: W^s(y, T_j) \cap T_k = \emptyset\}\;  \cup \;  \{
y \in V: W^s(y, T_j) \cap T_k \neq \emptyset\}$$ 1er. caso) Si el
primero de esos subconjuntos es vac\'{\i}o, entonces el otro es el
abierto $V$ y tiene entonces todos sus puntos interiores: $V
\subset A_{j,k}$.

2do caso) Si existe $y \in V$ tal que $W^s (y, T_j) \cap T_k =
\emptyset$ probemos que $y \in \mbox { int } \{ x \in V: W^s (y,
T_j) \cap T_k = \emptyset \} \subset A_{j,k}$:

En efecto, si as\'{\i} no fuera existir\'{\i}an $Y_n \rightarrow y$ en $T_j$
tales que $Z_n \in W^s(y_n, T_j) \cap T_k \neq \emptyset$. Tomando
subsucesiones convergentes, tendr\'{\i}amos $z_n \rightarrow z$ en $T_j
\cap T_k$ (porque $T_j, T_k$ son cerrados). Luego
$$z \in W^s_{\epsilon }(y_n) \Rightarrow \dist (f^p z_n, f^p y_n) \leq \epsilon
\; \; \forall p \geq 0$$ Por continuidad $\dist (f^pz, f^p y) \leq
\epsilon \; \forall p \geq 0$ y entonces $z \in W^s_{\epsilon }
(y)\cap T_j \cap T_k$. Luego $W^s(y,T_j) \cap T_k \neq \emptyset $
contra lo supuesto. $\; \; \Box$

 \vspace{.5cm}

{\bf \large 4.4 Demostraci\'{o}n de que ${\cal R}$ es una partici\'{o}n de
$M$}

\vspace{.3cm}

En la secci\'{o}n 4.2 se construy\'{o} una familia ${\cal R} =
\{\overline{R(x)} \;\}_{x \in Z^*}$ finita de rect\'{a}ngulos propios.
Para demostrar que ${\cal R} $ es una partici\'{o}n por cerrados de
$M$ falta probar que
\begin{itemize}
\item [ I)] ${\cal R}$ cubre $M$.
\item [II)] Los interiores de dos rect\'{a}ngulos distintos de ${\cal
R}$ son disjuntos.
\end{itemize}

\begin{proposition}
${\cal R } = \{ \overline {R(x)}\; \}_{x \in Z^*}$ definido seg\'{u}n
la secci\'{o}n 4.1 es un cubrimiento de $M$
\end{proposition}

{\em Demostraci\'{o}n: } Como la familia ${\cal R }$ es finita se
tiene:
$$\bigcup _{x \in Z^*} \overline {R(x)} \; = \overline{ \bigcup _{x \in Z^*} R(x)} \;
$$
Por la secci\'{o}n 4.2 si $x \in Z^*$ entonces $x \in R(x)$. Entonces
$\cup _{x \in Z^*} R (x) \supset Z^*$ de donde
$$\overline {\bigcup _{x \in Z^*} R(x)}\; \; \supset \; \overline{Z^*} \; \; = \; M $$
debido a la densidad de $Z^*$ en $M$. $\; \; \Box$

\begin{proposition}
Sea $R(x)$ definido en la secci\'{o}n 4.2. Se cumple
\begin{itemize} \em
\item [ 1) ] $\mbox { int } \overline{R(x)} \; \; = R(x)$
\item [ 2) ] \em
Si $x,y $ son tales que $R(x) \cap R(y) \neq \emptyset$ entonces
$R(x)= R(y)$.
\end{itemize}

\end{proposition}
{\em Demostraci\'{o}n : } 1) Por la definici\'{o}n en la secci\'{o}n 4.2 se
tiene:
$$R(x) = \bigcap _{(j,k,n) \in H(x) } \mbox { int } \overline{T_{j,k} ^n}\;$$
Por otro lado $R(x)$ es abierto contenido en $\overline{R(x)}\;$ y
entonces $R(x) \subset \mbox { int } \overline {R(x)}\;$.

Adem\'{a}s $R(x) \subset \cap_{(j,k,n) \in H(x)} \overline {T_{j,k}^n}
 \; $ que es un cerrado. Luego $\overline{R(x)}\; \subset \overline {T_{j,k}^n}
 \; \; \; \forall (j,k,n) \in H(x) $.

 Se deduce que
 $$\mbox{ int } \overline{R(x)} \; \subset \mbox { int }  \overline {T_{j,k}^n}
 \; \; \; \forall (j,k,n) \in H(x)$$ de donde

$$\mbox{ int } \overline{R(x)} \; \subset \bigcap _{(j,k,n) \in H(x)}
\mbox { int }  \overline {T_{j,k}^n}
 \; \; = R(x)$$

2) $R(x) cap R(y)$ es abierto no vac\'{\i}o y como $Z^*$ es denso en $M
$ contiene alg\'{u}n punto $z \in Z^*$. Alcanza probar que $R(x) =
R(z) = R(y)$. Para eso es suficiente demostrar que si $z \in Z^*$
con $z \in R(x)$ entonces $H(x) = H(z)$ (por la definici\'{o}n en la
secci\'{o}n 4.2).

Demostremos primero que $H(x) \subset H(z)$:

$(j,k,n) \in H(x) \Rightarrow x \in T_{j,k}^n \Rightarrow z \in
R(x) \subset \overline {T_{j,k}^n} \; \; \subset T_j =
\cap_{n=1}^4 T_{j,k}^n $. Luego existe $n_1$ tal que $z \in
T_{j,k}^{n_1}$. Pero $z \in R(x) = \overline {T_{j,k}^n} \; $,
luego
$$\mbox { int } T_{j,k}^{n_1} \cap \overline { T_{j,k} ^n} \; \neq \emptyset$$
de donde $n_1 = n$. Luego $z \in T_{j,k}^n$ y $(j,k,n) \in H(z)$
como se quer\'{\i}a probar.

Ahora probemos que $H(z) \subset H(x)$:

Sea $(i,h,m) \in H(z)$, lo que implica $z \in T_{i,h}^m \subset
T_i$, de donde $z \in T_i$. Sea $j $ tal que $x \in T_j = \cup
_{n=1}^4 T_{j,i}^n$. Existe $n$ tal que $x \in T_{j,i} ^n$. Como
$H(x) \subset H(z)$ entonces $z \in T_{j,i}^n \subset T_j$. Pero
entonces $z \in T_j \cap T_i = T_{j,i}^1$, o sea $n=1$. Luego $x
\in T_{j,i}^1 = T_j \cap T_i \subset T_i = \cup_{m=1}^4
T_{i,h}^m$.

Luego existe $m_1$ tal que $x \in T_{i,h} ^{m_1}$. Como $H(x)
\subset H(z) $ entonces $z \in T_{i,h}^{m_1}$. Pero por hip\'{o}tesis
$z \in T_{i,h}^m$ y entonces $m = m_1$ y $x \in T_{i,h}^m$ de
donde $(i,h,m) \in H(x)$ como se quer\'{\i}a probar. $\; \; \Box$

 \vspace{.5cm}

{\bf \large 4.5 Densidad de $Z^*$ en las variedades estable e
inestables}

\vspace{.3cm}

En la secci\'{o}n 4.3 se prob\'{o} que $Z^*$ es denso en $M$ y adem\'{a}s que
$$Z^* = \bigcap _{j,k = 1}^m T_j ^c \cup (A_{j,k} \cap B_{j,k})$$
Se probar\'{a} que $Z^*$ es adem\'{a}s denso en las $\epsilon$- variedades
estables e inestables que pasan por alg\'{u}n punto de $Z^*$. Este
resultado se utilizar\'{a} luego en la demostraci\'{o}n del teorema de
Sinai.
\begin{lemma}
Sean $A_{j,k}$ y $B_{j,k}$ definidos en la secci\'{o}n 4.3.

Si $x \in A_{j,k}$ entonces $\mbox{ int } T_j \cap W^s_{\epsilon}
(x) \subset A_{j,k}$.

Si $x \in B_{j,k}$ entonces $\mbox{ int } T_j \cap W^u_{\epsilon}
(x) \subset B_{j,k}$.
\end{lemma}

{\em Demostraci\'{o}n: } Si $x \in A_{j,k}$ existe un entorno $V$ de
$x$ en $T_j$ tal que para todo $y \in V$ $W^s(y, T_j)$ corta a
$T_k$ si y solo si lo hace $W^s(x, T_j)$ (por definici\'{o}n del
abierto $A_{j,k}$ en la secci\'{o}n 4.3).

Sea $w \in \mbox{ int } T_j \cap W^s_{\epsilon } (x)$. Tenemos que
$x = [w,x]$. Existe un entorno $W$ de $w$ en $M$ tal que $[W, x]
\subset V$ por la continuidad de la funci\'{o}n corchete. Podemos
elegir $W \subset T_j$ porque $w \in \mbox{ int } T_j$. Si $z \in
W $ entonces $[z,x] = y \in V$. Luego $W^s(z, T_j) = W^s(y, T_j)$.
Entonces  para todo $z \in W$: $W^s(z, T_j)$ corta a $T_k$ si y
solo si lo hace $W^s(x, T_j)$. De donde $w \in W \subset A_{j,k}$.

 An\'{a}logamente se tiene que $\mbox{ int } T_j \cap W^u_{\epsilon} (x) \subset B_{j,k}$
 cuando $x \in B_{j,k}$. $\;\; \Box$

 \begin{proposition}
 Sean $x \in Z^*$ y $R(x)$ definidos en la secci\'{o}n 4.2.

 $Z^*$ es denso en $W^s(x, R(x))$ y en $W^u(x, R(x))$.

 \end{proposition}

{\em Demostraci\'{o}n: } Dado un abierto $V$ en $W^s(x, R(x))$
 demostremos que contiene alg\'{u}n punto de $Z^*$. Si $w \in V$,
 entonces $w \in W^s_{\epsilon } (x) \cap R(x), \; w = [x,w]$. Por
 continuidad de la funci\'{o}n corchete existe un entorno $U$ de $w$
en $M$ tal que $[x,U] \subset V$.

Como $w \in R(x)$ y $R(x)$ es abierto podemos suponer que $U
\subset R(x)$.

Sabemos que $Z^*$ es denso en $M$ por lo demostrado en la secci\'{o}n
4.3. Existe $y \in U \cap Z^*$. Basta demostrar que el punto $z$
definido como $z = [x,y]$ pertenece a $Z^*$.

En efecto por lo demostrado en la secci\'{o}n 4.3: $Z^* = \cap
_{j,k=1} ^m T_j ^c \cup (A_{j,k} \cap B_{j,k})$. Basta demostrar
 que si $z \in T_j$ entonces $z \in A_{j,k} \cap B_{j,k} \; \; \forall
 k$.

 Por construcci\'{o}n $z = [x,y]\; \; , \; x,y \in Z^* \cap R(x)$. Sea
 $j $ tal que $z \in T_j$. Sea $i$ tal que $x \in T_i$. Se tiene
 que $y,z \in R(x) \subset \mbox{ int } \overline{ T_i} \; = \mbox { int }
  T_i$ pues $T_i $ es cerrado. As\'{\i} $x,y \in Z^* \cap T_i$. Por la
  caracterizaci\'{o}n de $Z^*$: $x,y \in A_{i,j} \cap B_{i,j}$.

  Adem\'{a}s $z \in W^s_{\epsilon }(x)\cap \mbox{ int } T_i, \; \; \;
  z \in W^u _{\epsilon } (y) \cap \mbox{ int } T_i$. Aplicando el
  lema de la secci\'{o}n 4.3 se tiene que $z \in A_{i,j}\cap B_ {i,j}
   = \cup _{n=1}^4 \mbox { int } T_{i,j} ^n$.

   Como $z \in T_j \cap T_i = T_{i,j}^1$ entonces $z \mbox{ int }
   T_{i,j}^1$. Como $x \in T_i$ entonces $x \in T_ {i,j}^n$ para
   alg\'{u}n $n$. Luego $y,z \in R(x) \subset overline{T_{i,j}^n} \;$.
   Luego $z \in \overline{T_{i.j}^n} \; \cap \mbox{ int } T_{i,j}
   ^1\neq \emptyset $ de donde $n = 1$.

   Hemos probado que para todo $j$ tal que $z \in T_j$ se cumple $x,y \in
   T_j$. Entonces $z \in T_j \Rightarrow x \in T_j \Rightarrow z \in R(x) \subset \mbox{ int }
    T_j$. Luego $z \in W^s_{\epsilon }(x) \cap \mbox { int } T_j \; \; \; , z \in W^u_{\epsilon }
    (y) \cap \mbox { int } T_j$. Como $x, y \in Z^* \cap T_j$
    entonces $x,y \in A_{j,k} \cap B_{j,k} \; \; \forall k$.
    Aplicando de nuevo el lema de la secci\'{o}n 4.3 se deduce $z \in A_{j,k} \cap B_{j,k} \;
    \; \; \forall k$. $ \; \; \Box$

    \begin{corollary}
    Si $x \in f^{-1} (Z^*)$ entonces $f^{-1} (Z^*)$ es denso en
    $W^s(x, R(x))$.
    \end{corollary}

{\em Demostraci\'{o}n: } Sea $V$ un abierto no vac\'{\i}o de $W^s(x,
R(x))$. Encontraremos un punto $y \in f^{-1} (Z^*) \cap V$.

 $V= U \cap W^s_{\epsilon } (x)$ donde $U$ es un abierto de $M$.

 $f(V) = f(U) \cap f(W^s_{\epsilon }(x))$ porque $f$ es
 invertible.

 $f(V) = f(U) \cap W^s_{\epsilon } (fx) \cap f(\overline {B}_{\epsilon
 }(x))$. Luego:

 $f(V)$ es un entorno no vac\'{\i}o en $W^s_{\epsilon }(fx)$, de donde
 tambi\'{e}n lo es $f(V) \cap R(fx)$ porque $R(fx)$ es un abierto de
 $M$.

 Por la proposici\'{o}n anterior existe $z \in Z^*$ en $f(V) \cap
 R(fx)$. Sea $y = f^{-1}(z)$. Por construcci\'{o}n $y \in f^{-1} (Z^*) \cap
 V$. $\; \; \Box$

\vspace{1cm}

\pagebreak

\begin{center}
{\Large Particiones de Markov para difeomorfismos de Anosov - }
Eleonora Catsigeras
\end{center}

\section{Teorema de Sinai}
En la parte 4 se construy\'{o} una partici\'{o}n ${\cal R}$ por
rect\'{a}ngulos propios. Porbaremos que ${\cal R}$ es una partici\'{o}n de
Markov de $M$ para el difeomorfismo $f$ de Anosovo con lo cual
quedar\'{a} probada la existencia de particiones de Markov.

\vspace{.5cm}

{\bf \large 5.1 Enunciado del teorema de Sinai}

\begin{theorem} \label{511}
Si $f$ es un difeomorfismo de Anosov en $M$ entonces, dado $\beta
>0$ existe una partici\'{o}n de Markov de $M$ para $f$ con di\'{a}metro
menor que $\beta$.
\end{theorem}

{\bf \large 5.2 Lemma}

\begin{lemma}
Sea $\tau = \{T_1, T_2, \ldots, T_m\}$ el cubrimiento definido en
la secci\'{o}n 3.5 con di\'{a}metro suficientemente peque\~{n}o y sean $Z^*,
{\cal R}$ definidos en la secci\'{o}n 4.2 a partir de $\tau$.

Si $x,y \in Z^*$ y si $y \in f^{-1}(Z^*) \cap R(x) \cap
W^s_{\epsilon }(x)$ entonces $fy \in R(fx)$.
\end{lemma}

{\em Demostraci\'{o}n: }  $R(f(x)) = \cap _{(j,k,n)} \in H(fx) \mbox {
int } \overline {T_{j,k}^n} \; $ por la definici\'{o}n en la secci\'{o}n
4.2.

 Basta probar que $fx \in T_{j,k}^n  \Rightarrow  fy \in \mbox { int } \overline {T_{j,k}^n}\;
 $. Probemos primero que $fx \in T_j \Rightarrow fy \in T_j$:

  $fx \in T_j \Rightarrow fx = \theta (\sigma q)$ con $q_1= p_j, \; q_o =
  p_h$. Luego $x = \theta (q) \Rightarrow x \in T_h$.

  $y \in W^s_{\epsilon }(x) \cap R(x) \subset W^s_{\epsilon }(x) \cap T_h = W^s(x,
  T_h)$.

  Aplicando la proposici\'{o}n de la secci\'{o}n 3.6 e obtiene $fy \in W^s(fx, T_j) \subset
  T_j$.

  Ahora probemos que $fx \in t_{j,k}^n \Rightarrow fy \in
  T_{j,k}^n$:

  $fx \in T_{j,k}^n \subset T_j \Rightarrow fy \in T_j = \cup _{n=1}^4
  T_{j,k}^{n}$. Por absurdo supongamos que $fy \in T_{j,k}^{n_1}$
  con $n_1 \neq n$. Por hip\'{o}tesis $y \in W_{\epsilon }^s (x)$.
  Aplicando  \ref{148} se tiene $fy \in W^s_{\epsilon } (fx)$.
  Luego de la proposici\'{o}n \ref{215} se obtiene $W^s (fx, T_j) = W^s (fy,
  T_j)$.
  La hip\'{o}tesis de absurdo y la definici\'{o}n de los rect\'{a}ngulos $T_{j,k}^n $
   en la secci\'{o}n 4.4 implican que estrictamente uno de los
   conjuntos $W^u(fx, T_j), \; W^u (fy, T_j)$ corta a $T_k$.

   Supongamos para fijar ideas que $W^u (fx, T_j) \cap T_k \neq \emptyset,
   \; \; W^s (fy, T_j) \cap T_k = \emptyset$.

Como $fx \in T_j$ por \ref{351} existe $q \in \Sigma (P)$ tal que
$(\sigma q )_0 = p_j, \; fx = \theta (\sigma q= f (\theta q))$.
Llamando $P_i = q_0$ se tiene $x = \theta (q)$. Luego $x \in T_i$.

Por hip\'{o}tesis $ y \in R(x)\cap W^S_{\epsilon }(x)$. De lo anterior
$R(x) \subset T_i$. Luego $y \in W^s(x, T_i)$.

Por lo supuesto, existe $fz \in W^u (fx, T_j) \cap T_k$. Aplicando
la proposici\'{o}n de la secci\'{o}n 3.6 se tiene que $z \in W^u (x,
T_i)$.

Como adem\'{a}s $fz \in T_k$, por \ref{351} existe $\overline q \in
\Sigma (P)$ tal que $\sigma (overline {q})_0 = p_k, \; fz = \theta
(\sigma \overline q) = f \theta \overline q$. Llamando $p_h =
\overline q_0$ se tiene $z \in \theta (\overline q)$. Luego $z \in
T_h$.

Por otro lado como $x,y \in T_i = \cup_{n=1}^4 T_{i,h}^n$, existen
$n_1, n_2$ tales que $x \in T_{i,h}^{n_1}, \; y \in
T_{i,h}^{n_2}$. Luego $y \in R(x) \subset
\overline{T_{i,h}}^{n_1}\;$. Como $y \in Z^*$ entonces $y \in
\mbox{ int } T_{i,h}^{n_2}$. Por lo observado en la secci\'{o}n 4.1
$n_1 = n_2$.

Sabiendo que $z \in W^u (x, T_i) \cap T_h \neq \emptyset$ y que
$x,y \in T_{i,h}^{n_1}$ se obtiene que existe $w' \in W^u (y, T_i)
\cap T_h \neq \emptyset$.

 Sea $w = [z, w'] \in T_i \cap T_h$ porque $z, w' \in T_i \cap
 T_h$. Como $w' \in W^u_{\epsilon } (y)$ se tiene $w= [z, w'] =  [z,y] \in T_i \cap
 T_h$.

 Si el di\'{a}metro de los rect\'{a}ngulos $T_i$ se elige suficientemente
 peque\~{n}o, por la continuidad uniforme de $f$ en $M$ compacta se
 cumple que $w \in W^u_{\epsilon }(y), \; \; w,y \in T_i \Rightarrow \dist (f^pW, f^py)
 \leq \epsilon \; \; \forall  p \geq 0, \; \; \dist (fw, fy) < \epsilon \Rightarrow
 fw \in W^u_{\epsilon}(fy)$.

 Como $w \in W^s (z, T_h)$ se tiene $f w \in W^s(fz, T_k)$.

 Luego $f w = [fz, fy] \in T_j$ porque $fy, fz \in T_j$, de donde
 $fw \in W^u (fy, T_j) \cap T_k \neq \emptyset$ contra lo
 supuesto.

 Se observa que se han utilizado hasta aqu\'{\i} todas las hip\'{o}tesis
 excepto que $y \in f^{-1}(Z^*)$ y que la misma demostraci\'{o}n puede
 realizarse permutando $x $ e $y$, pues todas las h\'{\i}p\'{o}tesis
 utilizadas son sim\'{e}tricas en $x$ e $y$. Entonces $y \in R(x)\cap Z^* \Rightarrow
 y \in R(y)\cap R(x) \neq \emptyset  \Rightarrow R(y) = R(x)$ y
 como $x \in Z^*$ entonces $x \in R(x) = R(y)$. Adem\'{a}s $y \in W^s_{\epsilon }(x)
 \Rightarrow x \in W^s_{\epsilon } (y)$ por definici\'{o}n. De all\'{\i}
 que la suposici\'{o}n del principio no era restrictiva.

 Finalmente demostremos que $fx \in T_{j,k}^n \Rightarrow fy \in \mbox{ int } \overline{ T_{j,k}}
  ^n$:

  Tenemos que $f x \in T_{j,k}^n \Rightarrow fy \in T_{j,k}^n $.
  Adem\'{a}s por hip\'{o}tesis $fy \in Z^*$ o sea $fy \in T_{j,k} ^n \Rightarrow
  fy \in \mbox{ int } T_{j,k} ^n \subset \mbox{ int } \overline {T_{j,k}}
  ^n\; \; \; \Box$.

\vspace{.5cm}

{\bf \large 5.3 Demostraci\'{o}n del teorema de Sinai}

\vspace{.3cm}

La partici\'{o}n ${\cal R}$ de $M$ construida en la secci\'{o}n 4.4 est\'{a}
formada por rect\'{a}ngulos propios y es un refinamiento del
cubrimiento $\tau$ construido en la secci\'{o}n 3.5. Luego:
$$\mbox{ diam } {\cal R} \leq \max_{T_i \in \tau } \mbox{ diam } T_i$$
Por la proposici\'{o}n \ref{351} $\mbox{ diam } T_i \leq 2 \beta$
donde $\beta$ es un n\'{u}mero positivo arbitrario.

Para terminar de demostrar el teorema de Sinai solo hace falta
verificar que ${\cal R}$ cumple las condiciones (i) y (ii) de la
definici\'{o}n de partici\'{o}n de Markov (\ref{242}), lo cual se
demuestra a continuaci\'{o}n:

\begin{proposition} \em
Si $x \in \mbox{ int } R_i \cap f^{-1} \mbox{ int } R_j$ con $R_i,
R_j \in {\cal R}$ entonces

(i) $fW^s(x, R_i) \subset W^s(fx, R_j)$

(ii) $fW^u(x, R_i) \supset W^u(fx, R_j)$
\end{proposition}

{\em Demostraci\'{o}n: } Por lo visto en \ref{148} $fW^s_{\epsilon
}(x) \subset W^s_{\epsilon } (fx)$. Para demostrar (i) alcanza
probar que $f W^s(x, R_i) \subset R_j$. Prob\'{e}moslo primero en el
caso particular que
$$x \in Z^* \cap f^{-1} Z^* \cap \mbox{ int } R_i \cap f^{-1} \mbox{ int } R_j$$

Por lo visto en la secci\'{o}n 4.4 $R(x ) = \mbox{ int } R(x), \; \;
R(fx) = \mbox{ int } R(fx)$. Por lo visto en la secci\'{o}n 4.5 $Z^*
\cap f^{-1} Z^*$ es denso en $W^s (x, R_i)$.  Dado $y \in W^s (x,
R_i)$ existe $y_n \rightarrow y $ con $y_n \in Z^* \cap f^{-1} Z^*
\cap W^s (x, R(x))$. Por el lema de la secci\'{o}n 5.2 se tiene $f Y_n
\in R(f(x)) \; \forall n$. Luego $f y  = \lim f  y_n \in \overline
{ R(fx)} = R_j$ o sea $f W^s (x, R_i) \subset R_j$ como quer\'{\i}amos
probar.

Ahora prob\'{e}moslo en general:

Si $x \in \mbox{ int } R_i \cap f^{-1} \mbox{ int } R_j$ sea
$$\overline x \in Z^* \cap f^{-1} Z^* \cap \mbox{ int } R_i \cap f^{-1} \mbox { int } R_j$$
Existe tal $x$ porque $Z^* \cap f^{-1} Z^*$ es denso en $M$ al ser
intersecci\'{o}n de abiertos densos.

$W^s(x, R_i) = \{ [x, \overline y]: \overline y \in W^s (\overline
x, R_i)\}\; \; \; f W^s(x, R_i) = \{ f[x, \overline y]: \overline
y \in W^s (\overline x, R_i)\} $. Como $f[x, \overline y] = [fx, f
\overline y]$ se obtiene:
$$fW^s(x, R_i)= \{ [fx, f \overline y]: \overline y \in W^s (\overline x, R_i)\}=
 \{[fx, w]: w \in fW^s(\overline x, R_i)\} $$
 $$ \subset \{[fx, w]: w \in W^s(f \overline x, R_j)\}= W^s(fx, R_j)$$

Hemos probado que $fW^s(x, R_i) \subset W^s (fx, R_j)$.

La afirmaci\'{o}n (ii) se prueba de (i) aplic\'{a}ndola al difeomorfismo
$f^{-1}$ recordando que las variedades estables de $f^{-1}$ son
las inestables de $f$. $ \; \; \; \Box$

 \vspace{15cm}

\pagebreak

\begin{center}
{\Large Particiones de Markov para difeomorfismos de Anosov - }
Eleonora Catsigeras
\end{center}

\section{Din\'{a}mica Simb\'{o}lica}

{\bf \large 6.1 Matriz de transici\'{o}n}

\vspace{.3cm}

Sea ${\cal R} = \{R_1, R_2, \ldots, R_m\}$ una partici\'{o}n de la
variedad $M$.

\begin{definition}
\em La matriz de transici\'{o}n $A$ de la partici\'{o}n es una matriz $m
\times m$ tal que $A_{i,j} = 1 $ sin $\mbox{ int } R_i \cap f^{-1}
\mbox{ int } R_j \neq \emptyset$ y $A_{i, j} = 0$ en caso
contrario.
\end{definition}

\begin{definition} \em
Si $A$ es la matriz de transici\'{o}n de la partici\'{o}n ${\cal R}$
denotaremos con $\Sigma _A$ al conjunto de sucesiones bi-infinitas
$\{R_{a_i}\}_{i \in \Z}$ de rect\'{a}ngulos de ${\cal R}$ tales que
dos rect\'{a}ngulos $R_{a_i}$ y $R_{a_{i+1}}$ consecutivos cumplen:
$$\mbox{ int } R_{a_i} \cap f^{-1} \mbox{ int } R_{\_{i+1}} \neq \emptyset$$
Luego
$$\Sigma _A = \{ a \in \{1,2,\ldots,m\}^{\Z}: A_{a_i a_{i+1}}=1\}$$
\end{definition}

\begin{remark} \em
\begin{itemize}
\item[ 1)] Si ${\cal R}$ es una partici\'{o}n de Markov, aplicando
\ref{261} se obtiene, para todo $x \in R_i \cap f^{-1} R_j$ cuando
$A_{ij} =1$:
\begin{itemize}
\item[ i)] $fW^s(x, R_i) \subset W^s(fx, R_j)$
\item[ii)] $fW^u(x, R_i) \supset W^u (fx, R_j)$
\end{itemize}
\item[ 2)] $\Sigma _A$ es invariante por el shift $\sigma$ pues
$$a \in \Sigma _A \Rightarrow A_{a_i a_{i+1}} =1 \; \forall i, \; \sigma (a)_i = a_{i+1}, \;
\sigma (a) _{i+1} = a_{i+2}, \; $$ $$A_{a_{i+1}a_{i+2}} =1 \;
\forall i \Rightarrow \sigma (a) \in \Sigma _A$$
\item[ 3)] En \ref{244} se observ\'{o} que si existe, es \'{u}nico el
punto $x \in \cap _{j \in \Z} f^{-j} R_{n_j}$ con $\{n_j\}_{j \in
\Z}$ sucesi\'{o}n cualquiera bi-infinita.
\end{itemize}
\end{remark}
Cuando ${\cal R} $ es una partici\'{o}n de Markov demostraremos que
\begin{itemize}
\item [ 1)] $a \in \Sigma _A \Rightarrow \; \exists x = \cap {j \in \Z} f^{-j} R_{a_j}$
\item [ 2)] La funci\'{o}n $\pi: \Sigma _A \mapsto M \; \;
\pi (x) = \cap _{j \in \Z} f^{-j} R_{a_j}$
es una semiconjugaci\'{o}n de $f$ con el shift.

$\pi$ llevar\'{a} entonces continuamente y sobreyectivamente las
\'{o}rbitas del shift en $\Sigma _A$ (llamada "din\'{a}mica simb\'{o}lica") en
las \'{o}rbitas de $f$ en $M$.
\end{itemize}
Adem\'{a}s en la secci\'{o}n 2.4 se observ\'{o} que si $x \not \in \cup _{j
\in Z} f^j \partial {\cal R}$ entonces existe una \'{u}nica sucesi\'{o}n
bi-infinita $\{n_j\}_{j \in \Z}$ tal que $f^j (x ) \in R_{n_j} \;
\forall j \in \Z$ (es decir $x = \cap _{j \in \Z} f^{-j}
R_{n_j}$). Eso significa que $\pi$ es adem\'{a}s inyectiva en los
puntos de $M \setminus \cup_{j \in \Z} f^j \partial {\cal R}$, que
como se vio en \ref{252} es denso en $M$.

 \vspace{.5cm}

{\bf \large 6.2 Lema}

\begin{lemma}
Sea ${\cal R} $ una partici\'{o}n de Markov y $A$ su matriz de
transici\'{o}n. Sea $a \in \Sigma _A$. Entonces:
$$K_N(a) = \bigcap_{j = -N}^{j = N} f^{-j} R_{a_j}$$ es un
rect\'{a}ngulo cerrado.
\end{lemma}

{\em Demostraci\'{o}n: }
 Supongamos en primer lugar que $K_N (a) \neq
\emptyset$ y demostremos la tesis en este caso: $K_N (a)$ es
cerrado por ser intersecci\'{o}n finita de cerrados $R_i$. Si $x,y \in
K_N(a) = \cap _{-N}^N f^{-j} R_{a_j}$ entonces $f^jx, f^j y \in
R_{a_j} \; \forall |j| \leq N$. Sea $w = [x,y] \subset R_{a_0}$.
Se cumple $w \in W^s(x, R_{a_0}), \; w \in W^u (y, R_{a_0})$.
Aplicando la proposici\'{o}n \ref{261} se obtiene:
$$fw \in W^s (fx, R_{a_1})\; \;, \; \; f^{-1} w \in W^u (f^{-1} y, R_{a_{-1}})$$
$$f^jw \in W^s (f^jx, R_{a_j})\; \;, \; \; f^{-j} w \in W^u (f^{-j} y, R_{a_{-j}})
\; \; \forall 0 \leq j \leq N$$ Se observa que para poder aplicar
la proposici\'{o}n \ref{261} se usan las hip\'{o}tesis $a \in \Sigma _A$ y
la partici\'{o}n ${\cal R}$ es de Markov.

Luego $ f^j w \in R_{a_j} \; \forall |j| \leq N$ o sea $w \in K_N
(a)$.

Demostremos ahora que $K_N(a) \neq \emptyset \; \; \forall a \in
\Sigma _A$:
$$K_N (a) = f^N R_{a_{-N}} \cap f^{N-1} R_{a_{-N +1}} \cap \ldots \cap f^{-N+1} R_{a_{N-1}}
 \cap f^{-N} R_{a_N} = $$ $$= f^N (R_{a_{-N}} \cap f^{-1} R_{a_{-N+1}} \cap
 \ldots \cap f^{-2N} R_{a_N})$$
Es suficiente demostrar que para todo $b \in \Sigma _A $ y para
todo $n \geq 0$
$$R_{b_0} \cap f^{-1} R_{b_1} \cap \ldots \cap f^{-n} R_{b_n} \neq \emptyset$$
Por inducci\'{o}n completa sobre $n$: Cuando $n= 0$ el conjunto
$R_{b_0} \neq \emptyset$ por ser un rect\'{a}ngulo.

Sabiendo por hip\'{o}tesis de inducci\'{o}n que existe $y \in R_{b_1} \cap
\ldots \cap f^{-n+1} R_{b_n}$ encontremos $z \in R_{b_0} \cap
f^{-1} (R_{b_1} \cap \ldots \cap f^{-n+1} R_{b_n})$:

$b \in \Sigma _A \Rightarrow \mbox{int} R_{b_0} \cap f^{-1}
\mbox{int } R_{b_1} \neq \emptyset \; \Rightarrow \; \exists
f^{-1} x \in R_{b_0} \cap f^{-1} R_{b_1}$.

Sea $z = [y,x]$. Se tiene $z \in R_{b_1}$ porque $x,y \in
R_{b_1}$. Adem\'{a}s $z \in W^s (y, R_{b_1})$. Por \ref{261}:
$$fz \in W^s (fy, R_{b_2}), \ldots, f^{n-1}z \in W^s (f^{n-1}y, R_{b_n})$$
As\'{\i} $z \in R_{b_1} \cap \ldots \cap f^{-n+1} R_{b_n}$. Adem\'{a}s $z
\in W^u (x: R_{b_1})$. Por \ref{261} $f^{-1}z \in W^u (f^{-1}x,
R_{b_0}) \subset R_{b_0}$. As\'{\i} se tiene
$$f^{-1}z \in R_{b_0} \cap f^{-1} (R_{b_1} \cap \ldots \cap f^{-n+1} R_{b_n})
\;\;\; \Box$$

\vspace{.5cm}

{\bf \large 6.3 Teorema de semiconjugaci\'{o}n}

\vspace{.3cm}

Sea ${\cal R} = \{R_i\}_{i= 1, \ldots, m }$ una partici\'{o}n de
Markov de la variedad $M$ para el difeomorfismo $f$.

Sea $A$ la matriz de transici\'{o}n y $\Sigma _A$ el subespacio de las
sucesiones bi-infinitas definido en la secci\'{o}n 6.1.

Entonces
\begin{itemize}
\item [ 1)] Para todo $a \in \Sigma _A$ existe \'{u}nico $x = \cap _{j \in \Z} f^{-j} (R_{a_j})$
\item [ 2)] La funci\'{o}n $\pi: \Sigma _A \mapsto M $ definida por $\pi (a) =
\cap _{j \in \Z} f^{-j} (R_{a_j})$ es una semiconjugaci\'{o}n de $f$
con el shift $\sigma$ en $\Sigma _A$.
\item [ 3)] $\pi |_{\pi^{-1} (\cap _{j \in Z} f^{-j} (\partial {\cal
R}) ^c)}$ es inyectiva, es decir $\pi$ es inyectiva en los puntos
de $M$ cuyas \'{o}rbitas no cortan al borde $\partial {\cal R}$ de la
partici\'{o}n.
\end{itemize}

{\em Demostraci\'{o}n: } Sea $K_N(a) = \bigcap _{|j| \leq N } f^{-j}
R_{a_j}$.

$K_N(a)$ es cerrado no vac\'{\i}o por el lema de la secci\'{o}n 6.2

$K_{N}(a) \supset K_{N+1}(a)$ por construcci\'{o}n.

Siendo $M$ compacta por la propiedad de las intersecciones finitas
se sabe que  el conjunto $$K(a) = \bigcap _{\-infty} ^{+\infty}
f^{-j} R_{a_j} = \bigcap _{N=1} ^{+\infty} K_N(a) \neq \emptyset$$
lo cual prueba que existe $x \in K(a)$.

$x \in K(a)$ es \'{u}nico porque si $x,y \in K(a)$ entonces $f^j (x),
f^j(y) \in R_{a_j} \; \forall j \in \Z$. Pero $R_{a_j}$ tiene
di\'{a}metro a lo sumo igual al de la partici\'{o}n de Markov que puede
elegirse menor que la constante de expansividad de $f$ ,
resultando $x = y$.

Se ha probado la parte 1) de la tesis.

Ahora probemos la parte 2):
$$\pi (\sigma (a)) = \bigcap _{j \in \Z} f^{-j} R_{a_{j+1}} = f (\bigcap _{j \in \Z} f^{-j-1} R_
{a_{j+1}}) = f(\pi (a))$$ Entonces $$\pi \circ \sigma (a) = f
\circ \pi (a) \; \; \forall \; a \in \Sigma _A$$

Por lo tanto es conmutativo el siguiente diagrama
$$\begin{array}{ccccc}
   & & \sigma &  &  \\
   & \Sigma _A & \mapsto & \Sigma _A &  \\
  \pi & \downarrow  &  & \downarrow  & \pi \\
   & M & \mapsto & M &  \\
    &   & f &   &   \\
\end{array}$$
Para demostrar que $\pi$ es una semiconjugaci\'{o}n hay que probar que
$\pi: \Sigma _A \mapsto M$ es continua y sobreyectiva.

$\pi$ es  continua pues si $a^n \rightarrow a \in \Sigma _A$,
llamando $\pi (a^n) = x_n \in M $ y eligiendo una subsucesi\'{o}n
convergente tenemos $x_n \rightarrow x \in M$ y adem\'{a}s:
$$x_n = \bigcap _{-\infty}^{+\infty} f^{-j} R_{a_j ^n}$$
Como $a^n \rightarrow a$, dado $p >0$ existe $N >0$ tal que $a_j^n
= a_j \; \forall n >N, \; \; \forall \; |j|\leq p$. Luego para
todo $n >N$ el punto $x_n \in \bigcap _{-p}^p f^{-j} R_{a_j}$ que
es un cerrado. Entonces $x = \lim x_n \in \bigcap _{-p}^p f^{-j}
R_{a_j}\; \; \forall \; p > 0$. As\'{\i} $x \in \bigcap _{-\infty }^{+
\infty} f^{-j} R_{a_j}= \pi (a)$.

$\pi $ es sobreyectiva pues si $x \in M \setminus \bigcup _{n \in
\Z} f^j \partial {\cal R}$ la \'{o}rbita por $x$ no corta al borde
$\partial {\cal R}$ de la partici\'{o}n. Entonces sea $a_j (x)$ el
\'{u}nico sub\'{\i}ndice tal que $f^j(x) \in \mbox{ int } R_{a_j(x)} $.
Como $f^j (x)\in \mbox{ int } R_{a_j(x)} \bigcap f^{-1} \mbox {
int } R_{a_{j+1} (x)}$ tenemos que $A_{a_j a_{j+1}} =1$ de donde
$a(x) \in \Sigma _A$.

Entonces por construcci\'{o}n: $$x \in \bigcap _{j \in \Z} f^{-j}
R_{a_j(x)} = \pi (a)\; \; \forall a \in \Sigma _A$$.

Hemos probado que $\pi (\Sigma _A) \supset M \setminus \bigcup_{n
\in Z} f^j \partial {\cal R}$. Por \ref{252} el conjunto $ M
\setminus \bigcup_{n \in Z} f^j \partial {\cal R}$ es denso en
$M$.

Adem\'{a}s $\Sigma _A$ es compacto porque es cerrado contenido en el
espacio compacto $\{1,2,\ldots, m\} ^{\Z}$. Luego $\pi (\Sigma
_A)$ es cerrado en $M$. Como contiene a un conjunto denso en $M$ y
es cerrado en $M$ es $M$, lo cual prueba que $\pi$ es
sobreyectiva.

Probemos ahora la parte 3): Si $a,a' \in \Sigma _A$ tales que $\pi
(a) = \pi (a') = x \in M \setminus \bigcup_{n \in Z} f^j \partial
{\cal R} $ entonces $f^j x \in \mbox{ int }R_{a_j} \cap \mbox{ int
} R_{a'_j} \; \; \forall \; j \in \Z$.

Por la definici\'{o}n de la secci\'{o}n 2.4 los rect\'{a}ngulos distintos
tienen interiores disjuntos. Entonces $a_j = a'_j$. Luego $a= a'$
y la transformaci\'{o}n $\pi$ restringida a la preimagen por $\pi $ de
$M \setminus \bigcup_{n \in Z} f^j \partial {\cal R}$, es
inyectiva. $\; \; \Box$

\vspace{.5cm}

{\bf \large 6.4 Conclusi\'{o}n}

\vspace{.3cm}

El teorema anterior permite construir una semiconjugaci\'{o}n $\pi$
del difeomorfismo de Ansosov $f$ con el shift $\sigma$ en el
subespacio de la sucesiones bi-infinitas $\Sigma _A $. Adem\'{a}s
$\pi$ es inyectiva en un conjunto denso en $M$.

Ya se observ\'{o} en la secci\'{o}n 2.4 que si ${\cal R }$ es un partici\'{o}n
por cerrados cualquiera de $M$, aunque no sea de Markov, existe
una funci\'{o}n $\pi$ sobreyectiva que puede demostrarse que es
continua usando la misma prueba de la parte 2) del teorema
anterior, tal que conmuta el siguiente diagrama:
$$\begin{array}{cccccc}
    &  & \sigma &   &   &   \\
    & \pi ^{-1} (M) & \mapsto & \pi ^{-1} (M) & \subset & \{1,2,\ldots,m\}^{\Z} \\
   & \pi \; \downarrow &   & \downarrow \; \pi &  &   \\
    & M  & \mapsto & M &   &   \\
    &   & f &   &   &   \\
\end{array}$$
Adem\'{a}s $\pi$ es inyectiva en $M \setminus \cup _{j \in \Z} f^{-j}
\partial {\cal R}$ que es denso en $M$.

En el caso que ${\cal R}$ sea adem\'{a}s una partici\'{o}n de Markov se
agrega que la semiconjugaci\'{o}n tiene dominio en $\Sigma _A$. El
subshift $\sigma |_{\Sigma _A} $ est\'{a} definido en el subconjunto
compacto de las sucesiones bi-infinitas que cumplen $A_{a_i
a_{i+1}} = 1$.

Se llama din\'{a}mica simb\'{o}lica a la din\'{a}mica del shift en $\Sigma
_A$. La existencia de una partici\'{o}n de Markov en $M$ para $f$
asegura la existencia de la din\'{a}mica simb\'{o}lica con la cual $f$ es
semiconjugada.

Finalmente se observa que en la definici\'{o}n de rect\'{a}ngulo, en la de
partici\'{o}n de Markov y en la demostraci\'{o}n del teorema de Sinai, no
se utiliza la diferenciabilidad de $f$ sino solo sus propiedades
topol\'{o}gicas. Es por lo tanto aplicable a una clase m\'{a}s general que
los difeomorfismos de Anosov: los homeomorfismos expansivos
topol\'{o}gicamente estables.


\begin{thebibliography}{2}



\bibitem[1]{1} R. Bowen: {\em Equilibrium states and the ergodic
theory of Anosov diffeomorphisms.} {\small Lecture Notes in Math.
{\bf 470.} Springer-Verlag } {1975}




\bibitem[2]{2} J. Lewowicz: {\em Lyapunov functions and Topological Stability.}
 {\small Journ. of Diff. Eq. { \bf 38}} {1980}

 \bibitem[3]{3} J. Lewowicz: {\em Invariant manifolds for regular points.}
 {\small Pacific Journ. of Math. { \bf 96}} {1981}

 \bibitem[4]{4} Hirsch-Pugh: {\em Stable manifolds and hyperbolic sets.}
 {\small Proc. Symp. in Pure Math. { \bf 14}} {1970}

 \bibitem[5]{5} R. Bowen: {\em Markov partitions for Axiom A diffeomorphisms.}
 {\small Amer. Journ. of Math.. { \bf 92}} {1970}

 \bibitem[6]{6} Y. Sinai: {\em Construction of Markov partitions.}
 {\small Funct. Anal. and its appl. { \bf 2}} {1968}

\end{thebibliography}
\end{document}